\documentclass[11pt,reqno]{amsart}
\usepackage{amsmath,amssymb,amsthm,amscd,mathrsfs,url,pinlabel,setspace,verbatim}
\usepackage[margin=1.25in,marginparwidth=1in,centering,letterpaper]{geometry}
\usepackage{xcolor,dcpic,latexsym,graphicx,epstopdf,comment}
\usepackage[all]{xy}
\usepackage{tikz,pgfplots,tikz-cd,circuitikz}
\usepackage{enumitem,upquote,tabularx,textcomp,makecell,caption}
\usepackage{extarrows,amsbsy,appendix,float,subfigure,fancyhdr,overpic,colonequals,microtype,diagbox}
\usepackage[numbers,sort&compress]{natbib}
\usepackage{hyperref}
\usepackage[all]{hypcap}

\xyoption{rotate}

\newtheorem {theorem}{Theorem}
\newtheorem {lemma}[theorem]{Lemma}
\newtheorem {proposition}[theorem]{Proposition}
\newtheorem {corollary}[theorem]{Corollary}
\newtheorem {assumption}[theorem]{Assumption}
\numberwithin{equation}{section}
\numberwithin{theorem}{section}

\theoremstyle{definition}
\newtheorem{definition}[theorem]{Definition}
\newtheorem{construction}[theorem]{Construction}
\newtheorem{notation}[theorem]{Notation}

\newtheorem{remark}[theorem]{Remark}

\newtheorem{warn}[theorem]{Warning}
\newtheorem{example}[theorem]{Example}
\newtheorem*{ack}{Acknowledgement}

\newcommand{\Z}{\mathbb{Z}}
\newcommand{\N}{\mathbb{N}}
\newcommand{\C}{\mathbb{C}}
\newcommand{\BN}{\mathrm{B}_\bullet\mathbb{N}}
\newcommand{\bX}{\mathbb{X}}
\newcommand{\rN}{\operatorname{N}_\bullet}

\newcommand{\BI}{\mathrm{BI}}
\newcommand{\bY}{\mathbb{Y}}
\newcommand{\bW}{\mathbb{W}}
\newcommand{\cA}{\mathcal{A}}
\newcommand{\cC}{\mathcal{C}}
\newcommand{\cD}{\mathcal{D}}
\newcommand{\sX}{\mathscr{X}}
\newcommand{\sD}{\mathscr{D}}
\newcommand{\cF}{\mathcal{F}}
\newcommand{\cO}{\mathcal{O}}
\newcommand{\cP}{\mathcal{P}}
\newcommand{\bK}{\mathbb{K}}
\newcommand{\fm}{\mathfrak{m}}
\newcommand{\ovv}[1]{{\overline{\overline{#1}}}}
\newcommand{\ov}[1]{{\overline{#1}}}
\newcommand{\un}[1]{{\underline{#1}}}
\newcommand{\wti}[1]{{\widetilde{#1}}}
\DeclareMathOperator*{\colim}{colim}
\DeclareMathOperator{\RG}{R\Gamma}
\DeclareMathOperator{\RGzd}{R\Gamma_{0d}}
\newcommand{\opp}{{\mathrm{opp}}}

\newcommand{\SHI}{\underline{\mathbf{SHI}}}

\newcommand{\inr}{\operatorname{int}}
\newcommand{\cone}{\operatorname{cone}}
\newcommand{\pt}{\mathrm{pt}}
\newcommand{\Ga}{\mathbf{G}_\mathrm{a}}

\newcommand{\Ch}{\mathsf{Ch}}
\newcommand{\proj}{\mathrm{proj}}
\newcommand{\Ndg}{\mathrm{N}^{\mathrm{dg}}}
\newcommand{\fin}{\mathrm{fin}}
\newcommand{\mix}{\mathrm{mix}}
\newcommand{\bbox}{\square}

\usetikzlibrary{calc,intersections}
\tikzset{every picture/.style=thick}
\pgfplotsset{compat=1.12}
\usetikzlibrary{arrows.meta,decorations.markings}

\newcommand{\bpf}{\begin{proof}}
\newcommand{\epf}{\end{proof}}
\newcommand{\bthm}{\begin{theorem}}
\newcommand{\ethm}{\end{theorem}}
\newcommand{\bprop}{\begin{proposition}}
\newcommand{\eprop}{\end{proposition}}
\newcommand{\bcor}{\begin{corollary}}
\newcommand{\ecor}{\end{corollary}}
\newcommand{\blem}{\begin{lemma}}
\newcommand{\elem}{\end{lemma}}
\newcommand{\bdefn}{\begin{definition}}
\newcommand{\edefn}{\end{definition}}
\newcommand{\bcons}{\begin{construction}}
\newcommand{\econs}{\end{construction}}
\newcommand{\bexmp}{\begin{example}}
\newcommand{\eexmp}{\end{example}}
\newcommand{\brem}{\begin{remark}}
\newcommand{\erem}{\end{remark}}
\newcommand{\benum}{\begin{enumerate}[leftmargin=*]}
\newcommand{\eenu}{\end{enumerate}}

\newcommand{\al}{\alpha}
\newcommand{\ga}{\gamma}
\newcommand{\la}{\lambda}
\newcommand{\de}{\delta}
\newcommand{\ep}{\epsilon}
\newcommand{\epp}{{\epsilon^{\prime}}}

\newcommand{\pa}{\partial}
\newcommand{\ot}{\otimes}%tensor product

\newcommand{\aand}{\ \text{and}\ }
\newcommand{\xra}{\xrightarrow}

\newcommand{\De}{\mathsf{D}}
\newcommand{\Mod}{\mathsf{Mod}}
\newcommand{\RMod}{\mathsf{RMod}}
\newcommand{\QCoh}{\mathsf{QCoh}}
\newcommand{\QCohzd}{\mathsf{QCoh}_{0\mathrm{d}}}
\newcommand{\K}{\mathsf{K}}
\newcommand{\Ac}{\mathsf{Ac}}
\newcommand{\Ani}{\mathsf{Ani}}
\newcommand{\Cat}{\mathsf{Cat}}
\newcommand{\Fun}{\mathsf{Fun}}
\newcommand{\End}{\mathsf{End}}
\newcommand{\Nil}{\mathsf{Nil}}

\DeclareMathOperator{\Eig}{Eig}
\DeclareMathOperator*{\limone}{lim^1}
\DeclareMathOperator{\gldim}{gl.dim}
\DeclareMathOperator{\Ext}{Ext}
\DeclareMathOperator{\Hom}{Hom}

\DeclareMathOperator{\id}{id}
\DeclareMathOperator{\fib}{fib}
\DeclareMathOperator{\cofib}{cofib}
\DeclareMathOperator{\Obj}{Obj}
\DeclareMathOperator{\Spec}{Spec}

\newcommand{\Pj}{\mathcal{P}}

\newcommand{\tbox}{\widetilde{\square}}

\title{Eigendecomposition for instantons via \texorpdfstring{$\infty$}{infinity}-category}
\author{Longke Tang}
\address{Simons Laufer Mathematical Sciences Institute}
\curraddr{}
\email{longket@alumni.princeton.edu}
\thanks{}

\author{Fan Ye}
\address{School of Mathematical Sciences\\Peking University}
\curraddr{}
\email{flyye@math.pku.edu.cn}
\thanks{}
\makeatletter

\begin{document}
\renewcommand\arraystretch{2}

\begin{abstract}
    Suppose $L$ is a framed link in a balanced sutured manifold $(M, \gamma)$. We construct a link spectral sequence in sutured instanton homology. For properly embedded surfaces in $M$ disjoint from $L$, we refine the spectral sequence to preserve the induced multi-grading. We use $\infty$-categories to encode homotopical information from instanton theory. In particular, we obtain a simultaneous eigendecomposition of homotopy-coherently commutative operators on a chain complex or a hypercube complex. As the algebraic foundation for subsequent work of the second author and collaborators, we construct the mixed hypercube by combining the hypercubes of positive and negative bypass maps. This construction is compatible with the eigendecomposition associated with the surface multi-grading on sutured instanton homology.
\end{abstract}

\maketitle
\vspace{-0.5\baselineskip}
\tableofcontents

\section{Introduction}
The surgery exact triangle in instanton Floer homology was first introduced by Floer \cite{floer1988instanton, floer1990knot} and later revisited by Braam--Donaldson \cite{donaldson95triangle} and Scaduto \cite{Scaduto2015instanton}. Analogous surgery triangles were subsequently established in monopole and Heegaard Floer homologies \cite{kronheimer2007monopoles, kronheimer2007monopolesandlens, Ozsvath2004c, ozsvath2005double,Freeman2021triangle,LRS2023triangle,AM2025sign,QY2026triangle}.

The basic setup consists of a framed oriented knot $K$ in a closed oriented smooth $3$-manifold $Y$, namely, a knot equipped with a chosen simple closed oriented curve $\lambda$ on the boundary $\partial N(K)$ of a tubular neighborhood of $K$. The curve $\lambda$, called the \emph{framing} of $K$, is chosen so that, for the oriented meridian $\mu=\pt\times\partial D^2\subset S^1\times D^2\cong N(K)$ of $K$, the pair $(\mu,\lambda)$ forms a basis of $H_1(\partial N(K);\mathbb Z)\cong\mathbb Z\oplus\mathbb Z$ and satisfies $\mu\cdot\lambda=1$. Equivalently, we have $\mu\cdot\lambda=-1$ when viewed on $Y\backslash \inr N(K)$.

For any Floer homology theory $F(-)$, there is a long exact sequence
\begin{equation}\label{eq: long exact sq}
    \cdots \to F(Y)\to F(Y_0(K))\to F(Y_1(K))\to F(Y)\to\cdots
\end{equation}
where $Y_i(K)$, for $i=0,1$, denotes the closed $3$-manifold obtained from $Y$ by Dehn surgery with coefficient $i$. Namely, $Y_i(K)$ is obtained by gluing a solid torus to $Y\backslash \inr N(K)$ so that the meridian $\pt\times\partial D^2$ of the solid torus is identified with the curve representing the homology class $i\mu+\lambda$. The maps in the exact sequence are induced by $4$-dimensional cobordisms given by attaching a single $2$-handle on the product of the source $3$-manifold and the interval $I$, possibly together with additional data arising from local coefficients or principal bundles. These maps are usually referred to as \emph{surgery cobordism maps}.

In fact, a stronger statement holds at the chain level. For the chain complex $CF(Y)$ defining $F(Y)$, one has a quasi-isomorphism in the derived category
\begin{equation}\label{eq: cone}
    CF(Y)\cong \cone\big(CF(Y_0(K))\to CF(Y_1(K))\big),
\end{equation}where $\cone$ denotes the mapping cone complex.

Since the long exact sequence \eqref{eq: long exact sq} is $3$-periodic, it can be viewed as an exact triangle. The quasi-isomorphism \eqref{eq: cone} is likewise $3$-periodic. We emphasize that the mapping cone formulation implies the existence of the exact triangle, whereas the converse implication generally requires working over a field.

The mapping cone description is particularly well suited for generalizations to surgeries on framed links. Let $L=\bigcup_{i=1}^nK_i$ be a framed link in $Y$. For $\epsilon\in\{0,1\}^n$, let $Y_\epsilon(L)$ denote the manifold obtained by performing surgery on $L$, with surgery coefficient $\epsilon_i$ on the component $K_i$. Then \eqref{eq: cone} generalizes to the quasi-isomorphism
\begin{equation}\label{eq: iterated cone}
    CF(Y)\cong \big(\bigoplus_{\epsilon\in\{0,1\}^n}CF(Y_\epsilon(L)),\sum_{\epsilon\leq\epsilon'\in\{0,1\}^n}D_{\epsilon,\epsilon'}\big),
\end{equation}
where the partial order is induced by $0<1$. Here we suppress the homological grading shifts on $CF(Y_\epsilon(L))$ for simplicity. The filtration associated to the $L^1$-norm $|\epsilon|$, which counts the number of entries equal to $1$ in $\epsilon$, gives rise to a spectral sequence
\begin{equation}\label{eq: spectral sequence}
    E_1=\bigoplus_{\epsilon\in\{0,1\}^n}F(Y_\epsilon(L))\Longrightarrow F(Y).
\end{equation}We will call such a construction a \emph{spectral sequence from the cube filtration}.

The link spectral sequence arising from this cube filtration was first established by Ozsv\'{a}th--Szab\'{o} \cite{ozsvath2005double} in Heegaard Floer theory over $\mathbb Z/2$ coefficients. Analogous constructions for monopole and instanton Floer homologies were later developed by Bloom \cite{Bloom2011link} over $\mathbb Z/2$ coefficients and by Scaduto \cite{Scaduto2015instanton} over $\mathbb Z$ coefficients, respectively. Recently, the monopole version was refined to $\mathbb Z$ coefficients by Qiu and the second author \cite{QY2026triangle}. 

A typical application is a spectral sequence from (even or odd, reduced or unreduced) Khovanov homology $Kh(L)$ of a link $L\subset S^3$ to an appropriate Floer homology $F(-\Sigma_2(L))$ of the double branched cover of $L$ with reversed orientation. We also refer the reader to \cite{kronheimer2011khovanov,kronheimer2021barnatan,francesco2019spectral,invol2023spectral,jiakai2024spectral,bhat2023newtriangle} for related constructions in other settings.

In this paper, we construct a variation of the above link spectral sequence for balanced sutured manifolds $(M,\gamma)$ (cf.\ Definition \ref{defn: balanced sutured manifold}). Here $(M,\gamma)$ consists of a compact oriented $3$-manifold $M$ with nonempty boundary and a collection of oriented curves $\gamma\subset\partial M$, called the \emph{suture}, satisfying certain balanced conditions. Kronheimer--Mrowka \cite{kronheimer2010knots} introduced two Floer homology theories, sutured monopole homology $SHM(M,\gamma)$ and sutured instanton homology $SHI(M,\gamma)$, as analogues of Juh\'{a}sz's sutured Floer homology $SFH(M,\gamma)$ in Heegaard Floer theory.

Both sutured homology theories are constructed by embedding $M$ into a closed oriented $3$-manifold $Y$ determined by $\gamma$ and certain auxiliary choices, such that the components of $\partial M\backslash\gamma$ lie on two parallel copies of an embedded oriented surface $R\subset Y$. One then extracts a subspace of $F(Y)$ associated to $R$ in the corresponding Floer theory. The pair $(Y,R)$ is called a \emph{closure} of $(M,\gamma)$. The resulting homology groups are independent, up to isomorphism, of the choice of closure. Further naturality results in the projective setting were established by Baldwin--Sivek \cite{BS2015naturality}, and the functoriality results were studied by Baldwin--Sivek \cite{baldwin2016contact,BS2016instanton} and Li \cite{li2018contact,li2018gluing}.

Moreover, given properly embedded surfaces $S_1,\ldots,S_k\subset M$ satisfying certain admissibility conditions (cf.\ Definition \ref{defn: admissible surface}), Ghosh--Li \cite{li2019decomposition} constructed a $\mathbb Z^k$-grading on $SHM$ and $SHI$, which we refer to as \emph{surface grading}. See also \cite{kronheimer2010instanton, li2019direct} for related constructions.

More precisely, in the construction of $SHM$, the subspace associated to the surface $R$ is obtained by evaluating the first Chern class of the spin$^c$ structure of $Y$ on $R$ and taking the summand corresponding to the maximal nontrivial evaluation, namely $2g(R)-2$. The surface grading is constructed similarly by considering evaluations on closed oriented surfaces $\overline{S}_l\subset Y$ containing $S_l$ as subsurfaces for $l=1,\ldots,k$. The analogous construction for $SHI$ is realized by taking simultaneous generalized eigenspaces of the operators $\mu(R)$ and $\mu(\overline{S}_l)$ acting on the instanton Floer homology of $Y$ with $\C$ coefficients. Hence $SHI$ is only defined over $\C$. In both cases, the possible evaluations or the possible eigenvalues are even integers with absolute values bounded above by $2g(\ov{S}_l)-2$ when $\ov{S}_l$ is connected. 

The main goal of this paper is to combine the link surgery spectral sequence \eqref{eq: spectral sequence}, sutured instanton homology $SHI$, and the surface grading structure. More precisely, we prove the following main theorem.

\bthm\label{thm: main grading for cube}
    Suppose $(M,\ga)$ is a balanced sutured manifold, and let $S_l\subset (M,\ga)$, $l=1,\ldots,k$, be admissible surfaces. Suppose that $L=\bigcup_{i=1}^nK_i\subset M$ is a framed link disjoint from all the surfaces $S_l$. For $\epsilon\in\{0,1\}^n$, let $M_\epsilon(L)$ be the sutured manifold obtained from $(M,\ga)$ by performing surgeries on $L$ with coefficients given by $\epsilon$. We also use $\ga$ and $S_l$ to denote the induced sutures and induced surfaces in $M_\epsilon(L)$, respectively. Then there exists a spectral sequence from the cube filtration
    \begin{equation*}\label{eq: spectral sequence SHI}
        \bigoplus_{\epsilon\in\{0,1\}^n}SHI(M_\epsilon(L),\ga)\Longrightarrow SHI(M,\ga)_L,
    \end{equation*}
    which respects the surface multi-grading associated to the surfaces $S_l$. Here the subscript $L$ indicates the additional bundle set used in the construction of instanton Floer homology (cf.\ the third line of \cite[Figure 1]{Scaduto2015instanton}). Moreover, the differential on the first page is given by the restrictions of the surgery cobordism maps on the generalized eigenspaces.
\ethm

% \cred{Actually, the spectral sequence in Theorem 1.1 comes from the fact that $SHI(M)$ is the total cofiber of an $n$-cube consisting of $SHI(M_\epsilon(L),\gamma)$. Add this to the introduction.}

\brem
    The analogous result of Theorem \ref{thm: main grading for cube} for $SHM$ follows naturally from the construction of the spectral sequence together with the surface grading. Here we consider $SHM$ without the additional bundle data $L$ over $\mathbb{Z}/2$ coefficients \cite{Bloom2011link}, or with the additional local system data $L$ over $\mathbb{Z}$ coefficients \cite{QY2026triangle}. Indeed, the evaluation of the spin$^c$ structures is compatible with all differentials appearing in the construction of the spectral sequence. On the other hand, the existence of the spectral sequence in Theorem \ref{thm: main grading for cube} for $n = 2$  without considering the surface grading follows from Ghosh--Zemke \cite[\S 6]{GZ2021}. However, due to the nontrivial intersections between the closed embedded surfaces $R$ and $\overline{S}_l$ arising from their constructions, verifying the compatibility of the spectral sequence with the surface grading associated to $S_l$ requires additional considerations beyond those in \cite[Remark 2.2 and \S 8.1]{GZ2021}. In particular, the results in \cite[\S 9]{GZ2021} can be viewed in the context of Theorem \ref{thm: main grading for cube}, together with the construction of homotopies and higher homotopies between $\mu$-operators in \cite[\S 8]{Ye2026infty}.
\erem
\brem\label{rem: dependence of data}
A priori, the spectral sequence in Theorem \ref{thm: main grading for cube} depends on choices arising from two sources:
\begin{enumerate}
    \item The choice of closure of $(M,\ga)$. The dependence on this choice may be studied using the naturality results of Baldwin--Sivek \cite{BS2015naturality}.
    \item The data used to construct the $\mu$-operators, as well as the homotopies and higher homotopies among them, as discussed in \cite[\S 8]{Ye2026infty}. The dependence on these choices requires additional care.
\end{enumerate}
\erem

The proof of Baldwin--Sivek \cite[\S 4]{BS2022khovanov} shows that the bypass exact triangle is a special case of the surgery exact triangle; a minor technical correction to the proof is provided in \cite[\S 8]{LY2025dimension}. Therefore, by modifying the statement of Theorem \ref{thm: main grading for cube}, we obtain the following theorem.

\bthm\label{thm: main bypass cube}
    Suppose $(M,\ga)$ is a balanced sutured manifold, and let $S_l\subset (M,\ga)$, $l=1,\ldots,k$, be admissible surfaces. Suppose that $\eta=\bigcup_{i=1}^n\eta_i\subset\partial M$ is a collection of pairwise disjoint arcs, called \emph{bypass arcs}, such that for any $i\in[1,n]$ and $j\in[1,k]$, we have
    \[
    \eta_i\cap S_j=\varnothing,\qquad |\eta_i\cap\ga|=3,
    \]
    and $N(\eta_i)\subset\partial M$ is as shown in the top subfigure of Figure \ref{fig:bypass}. For $\epsilon\in\{0,1\}^n$, let $\ga_\epsilon(\eta)$ be the suture obtained from $\ga$ by replacing the suture in $N(\eta_i)$ with one of the two sutures shown in the bottom two subfigures of Figure \ref{fig:bypass}, according to the value of $\epsilon_i$. Then there exists a spectral sequence from the cube filtration
    \begin{equation*}\label{eq: spectral sequence SHI bypass}
        \bigoplus_{\epsilon\in\{0,1\}^n}SHI(M,\ga_\epsilon(\eta))\Longrightarrow SHI(M,\ga)
    \end{equation*}
    that respects the surface multi-grading associated to the surfaces $S_l$. Moreover, the differential on the first page is given by the bypass maps.
\ethm
\newcommand{\BypassChords}[3]{%
  \begin{scope}[shift={(#1,#2)},rotate=#3]
    \foreach \a/\b in {0/180,60/120,-60/-120}{%
      \draw (\a:\BypassRadius) -- (\b:\BypassRadius);
    }
  \end{scope}%
}

\begin{figure}[htbp]
\centering
\resizebox{0.35\textwidth}{!}{%
\begin{tikzpicture}[line width=0.8pt,line cap=butt,line join=round,every node/.style={font=\LARGE}]
\def\BypassRadius{1.75}

\draw[dashed] (11,14) circle[radius=\BypassRadius];
\draw[dashed] (8.25,9.75) circle[radius=\BypassRadius];
\draw[dashed] (13.75,9.75) circle[radius=\BypassRadius];
\draw[dashed] ({11-\BypassRadius*cos(30)},14) -- ({11+\BypassRadius*cos(30)},14);

\BypassChords{11}{14}{90}
\BypassChords{8.25}{9.75}{-30}
\BypassChords{13.75}{9.75}{30}

\node at (11.5,13.5) {$\eta_i$};
\node at (11,11.5) {$\ga$};
\node at (8,7.25) {$\ga_0$};
\node at (13.75,7.25) {$\ga_1$};
\end{tikzpicture}%
}
\caption{The neighborhood of the dashed bypass arc $\eta_i$ and the two modifications $\ga_0$ and $\ga_1$ of the suture. We assume the orientation of $\partial M$ points out of the page.}
\label{fig:bypass}
\end{figure}
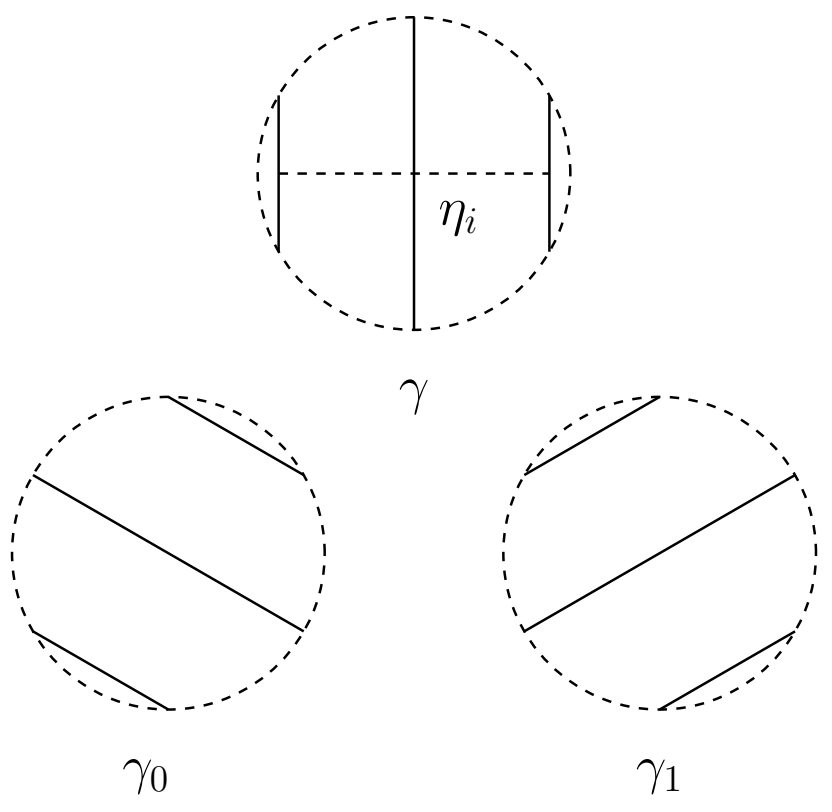
\brem
    The cyclic order $\ga,\ga_0,\ga_1$ in Figure \ref{fig:bypass} is opposite to the original convention in \cite[Figure 11]{BS2022khovanov}, since no orientation reversal is performed on $(M,\ga)$ (cf.\ Term (\ref{g2}) of Proposition \ref{prop: grading properties SHI}).
\erem

\brem\label{rem: further direction}
    Beyond Theorems \ref{thm: main grading for cube} and \ref{thm: main bypass cube}, another motivation of this paper is to develop the foundational constructions for future work by the second author and collaborators; see \S \ref{subsec: combination of bypass maps} for a sketch explanation of the construction. In that work, the authors plan to extend the ideas of \cite[\S 4]{BLY2020} to the setting of collections of disjoint bypass arcs. The goal is to construct a spectral sequence from the chain complex $SFC(M,\gamma)$ of $SFH(M,\gamma)$ to $SHI(M,\gamma)$ by combining both positive and negative bypass maps, such that the differentials affect the surface grading in a controlled way.
\erem
\begin{ack}
    We would like to thank Ian Zemke for helpful conversations on eigendecompositions. The second author thanks Michael Hopkins, David Treumann, Leon Liu, and Yuanning Zhang for helpful conversations on $\infty$-categories. The second author would also like to thank John A. Baldwin, Zhenkun Li, and Steven Sivek for discussions related to upcoming collaborative work, which provided motivation for this project. The second author was partially supported by Simons Collaboration Grant \#271133 through Peter Kronheimer during the preparation of this work.
\end{ack}

\noindent\textbf{AI Declaration.}
An early version of this work was completed in April 2025. The authors used OpenAI Codex with the \texttt{gpt-6-astra} and \texttt{gpt-6.1-sol} models to correct typos in the text and mathematical notation and to polish the presentation of some statements. The authors reviewed the resulting changes and take full responsibility for the manuscript.

\section{Sketch of the constructions}
This section is a continuation of the introduction. We provide a sketch of the constructions of hypercubes and the eigendecomposition, explain why we use $\infty$-categories, and describe the construction combining both positive and negative bypass maps. Finally, we provide an outline of the remaining sections. 

% The following is the notation used in the rest of the paper.
% \benum
%     \item For $n\in\Z_+$, we will write $\un{n}=\{1,\ldots,n\}$.
%     \item We usually use indices $i,j\in\un{n}$, $l\in\un{k}$, and $p,q\in\un{r}$.
%     \item Let $\bK$ denote the ground commutative ring with identity. It could be $\Z,\Q,\C$, and $\C[\{U_n\}_{i\in\un{n}}]$.
%     \item For $i\in\N$, we write $\vec{i}=(i,\ldots,i)$.
%     \item If not mentioned, we choose the interval $I$ to be $[0,1]$.
% \eenu
\subsection{Hypercube and eigendecomposition}\label{subsec: Algebraic setup of spectral sequences}
The construction on the right-hand side of \eqref{eq: iterated cone} gives a \emph{hypercube complex (of dimension $n$)}, or simply a \emph{hypercube}. The filtration on this hypercube induced by the $L^1$-norm $|\cdot|$ of the indices gives rise to a spectral sequence. This filtration is commonly referred to as the \emph{cube filtration}. See \cite[\S 5]{MO2022link}, \cite[\S 2.2]{Zemke23bordered}, and \cite{BHL2019cancel}. More precisely, a hypercube is a chain complex of the form
\begin{equation}\label{eq: hypercube}
(C,D)=\left(\bigoplus_{\ep\in\{0,1\}^n}C_\ep[|\ep|-n],~\sum_{\ep\leq\epp}D_{\ep,\epp}\right),
\end{equation}
where $C_\ep$ are finitely presented projective modules over a commutative Noetherian ring $\bK$ of finite global dimension, $D_{\ep,\epp}:C_\ep\to C_\epp$ are module morphisms, and the partial order on $\{0,1\}^n$ is induced by $0<1$. Each $C_\ep$ is $\Z$- or $\Z/m$-graded for some $m\in\Z_+$, and $D_{\ep,\epp}$ has degree $|\epp-\ep|-1$ before shifting. We use $(C_\ep[k])_j=(C_\ep)_{j+k}$, so $D$ has total degree $-1$. For periodic gradings assume that $m$ is even or $2=0$ in $\bK$.

We use the total differential signs of \cite[Construction 7.7]{Ye2026infty}, with $D_{\ep,\ep}=(-1)^{n-|\ep|}\partial_\ep$ for the differential $\partial_\ep$ on the unshifted vertex complex. The compatibility condition is
\begin{equation}\label{eq: compatibility}
\sum_{\ep\leq\de\leq\epp}D_{\de,\epp}\circ D_{\ep,\de}=0,
\end{equation}
or equivalently $D\circ D=0$, which implies that $(C,D)$ is a chain complex. In the Floer setting, $(C_\ep,\partial_\ep)$ is the chain complex defining $F(Y_\ep(L))$ in \eqref{eq: spectral sequence}, where the homological shifts are suppressed. The other morphisms $D_{\ep,\epp}$ for $\ep<\epp$ in \eqref{eq: iterated cone} are induced by families of cobordisms and are called \emph{family cobordism maps} in the Floer homology setting.

An operator
\[
\mu=\sum_{\ep\leq\epp\in\{0,1\}^n}\mu_{\ep,\epp}:(C,D)\to(C,D)
\]
on a hypercube $(C,D)$ consists of module morphisms $\mu_{\ep,\epp}:C_\ep\to C_\epp$ satisfying $\mu\circ D=D\circ\mu$, with each component of degree $|\epp-\ep|$ on the unshifted vertices, so $\mu$ has total degree zero.

If $\bK=\mathbb C$ (or more generally an algebraically closed field), then linear algebra gives the generalized eigendecomposition
\begin{equation}\label{eq: eigenspace decomp}
(C,D)=\bigoplus_{\la\in\bK}(C^\la,D^\la),
\end{equation}
where
\[
C^\la=\colim_{N\to\infty}\ker(\mu-\la\cdot\id)^N,\qquad D^\la=D|_{C^\la}.
\]
However, the subcomplex $(C^\la,D^\la)$ is generally not a hypercube complex, but only a $\{0,1\}^n$-filtered chain complex. In particular, $C^\la$ need not be the direct sum of the shifted modules $C_\ep^\la[|\ep|-n]$, where
\[
C_\ep^\la=\colim_{N\to\infty}\ker(\mu_{\ep,\ep}-\la\cdot\id)^N
\]
for $\ep\in\{0,1\}^n$.
\begin{example}
    When $n=1$, the total complex has underlying graded module $C_0[-1]\oplus C_1$ and is the mapping cone
\[
\cone(D_{0,1}:(C_0,\partial_0)\to(C_1,\partial_1)).
\]
An element in $C^\la$ is either $(0,x_1)\in0\oplus C_1^\la$ or $(x_0,x_1')\in C_0^\la[-1]\oplus C_1$, but $(x_0,0)$ need not belong to $C^\la$, since $\mu_{0,1}$ may be nonzero.
\end{example}

To construct a hypercube structure on the eigendecomposition, Ghosh--Zemke \cite[\S 6]{GZ2021} show that the associated graded module of the $\{0,1\}^n$-filtered module $C^\la$ is canonically isomorphic to \[\bigoplus_{\ep\in\{0,1\}^n}C_\ep^\la[|\ep|-n].\] Given a choice of an isomorphism called a \emph{splitting map}
\begin{equation}\label{eq: map phi}
\Phi:\bigoplus_{\ep\in\{0,1\}^n}C_\ep^\la[|\ep|-n]\to C^\la
\end{equation}
inducing the inverse of this canonical identification, one can construct a hypercube complex
\[
E^\la(C,\mu;\Phi)=\left(\bigoplus_{\ep\in\{0,1\}^n}C_\ep^\la[|\ep|-n],~\Phi^{-1}\circ D\circ\Phi\right).
\]

Moreover, in \cite[\S 6]{GZ2021}, this construction was extended to morphisms between hypercubes and applied to the simultaneous generalized eigenspaces of two homotopy-coherently commutative operators $\mu_1,\mu_2$ on $(C,D)$. Here homotopy-coherent commutation means that there exists a morphism
\[
h=\sum_{\ep\leq\epp\in\{0,1\}^n}h_{\ep,\epp}:C\to C,
\]
often called a \emph{homotopy}, such that the following diagram
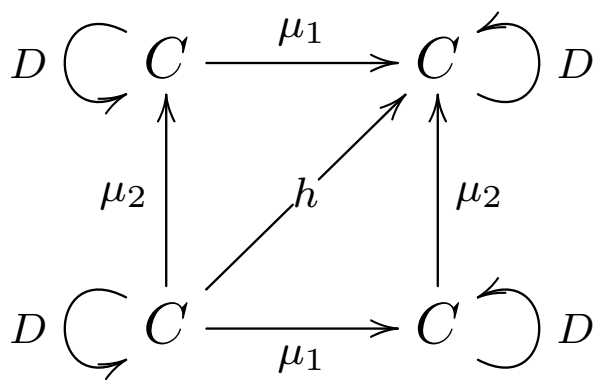
\begin{equation}\label{eq: homotopy commute operator 2}
\xymatrix@C=3pc@R=3pc{C\ar@(lu,ld)_{D}\ar[r]^{\mu_1}&C\ar@(rd,ru)_{D}\\
    C\ar@(lu,ld)_{D}\ar[r]_{\mu_1}\ar[u]^{\mu_2}\ar[ur]|-{h}&C\ar@(rd,ru)_{D}\ar[u]_{\mu_2}}
\end{equation}is homotopy commutative. On the total complex $(C,D)$, the degree-$1$ homotopy $h$ satisfies
\[
\mu_2\mu_1-\mu_1\mu_2=Dh+hD.
\]
With the signs of \cite[Construction 7.7]{Ye2026infty}, these data define an $(n+2)$-dimensional hypercube satisfying \eqref{eq: compatibility}.

While it is possible to extend this strategy to more homotopy-coherently commutative operators $\mu_1,\ldots,\mu_r$, it becomes increasingly subtle to choose the isomorphisms $\Phi$ in \eqref{eq: map phi} consistently and to keep track of the associated homotopies and higher homotopies. In this paper, instead of carrying out this process explicitly, we use $\infty$-categories to systematically organize the higher homotopy-coherent information. 

From the perspective of $\infty$-categories, a hypercube is simply a diagram indexed by the category $\bbox^n=(\Delta^1)^n$ (see Definition \ref{defn: intro}), and homotopy-coherently commutative operators correspond to actions of the classifying $\infty$-category $\BN^r$ associated to the monoid $\N^r$ on such a diagram. See \cite[\S 2.1--2.2]{Ye2026infty} for the corresponding constructions in instanton theory. The equivalence between these two perspectives is established in Theorem \ref{thm: main infty replacement}.

\subsection{Reason for using \texorpdfstring{$\infty$}{∞}-categories}\label{subsec: language infty cate}
To account for homotopies and higher homotopies, we will extensively employ the theory of ($\infty$,$1$)-categories as developed by Lurie \cite{HTTlurie2009, HAlurie2017, kerodon}. This framework allows us to systematically capture the information encoded by hypercubes and homotopy-coherently commutative operators and to develop the eigendecomposition. In line with Lurie's terminology, we use the term ``$\infty$-categories'' to refer to ($\infty$,$1$)-categories. For readers who are unfamiliar with $\infty$-categories, standard references include \cite[Chapter 1]{HTTlurie2009}, \cite[Chapter 1]{HAlurie2017}, \cite[\S 1-2]{lurie2009cob}, and the survey articles \cite{Bergner2010survey1,Bergner2020survey2}.

However, unlike in those references (and also unlike the previous work of the second author \cite{Ye2026infty}), we will use the language of $\infty$-categories in a model-independent manner; that is, we will not identify $\infty$-categories with any specific model, such as quasi-categories, complete Segal spaces, or simplicial categories. This approach stems from our focus not on \emph{developing} the theory of $\infty$-categories but on \emph{utilizing} it. The model-independent style aligns more closely with our intuition: we consider $\infty$-categories as categories while remembering all the higher coherence data crucial to their application. The models merely \emph{represent} $\infty$-categories, akin to how Dedekind cuts and Cauchy sequences \emph{represent} real numbers; we do not typically conceive of real numbers as Dedekind cuts or Cauchy sequences when working with them. This approach is consistent with the writing style of many contemporary papers utilizing $\infty$-categories, such as \cite{prisms2022, BurklundHahnSenger2023}.

To distinguish an object of the $\infty$-category of spaces from an actual topological space, we often refer to the former as an \emph{anima}, and use $\Ani$ to denote the $\infty$-category of animas (\cite[Definition 1.2.16.1]{HTTlurie2009}). We use $\Cat$ to denote the $\infty$-category of small $\infty$-categories (\cite[Definition 3.0.0.1]{HTTlurie2009}), and $\Fun$ to denote the $\infty$-category of functors (\cite[Notation 1.2.7.2]{HTTlurie2009}). Generally, we prefer to use sans serif letters to denote fixed $\infty$-categories or constructions that output $\infty$-categories, while calligraphic letters denote general $\infty$-categories.

%Sometimes, we also use a \emph{$k$-simplex} to refer to a $k$-morphism and use a \emph{simplex} to refer to either a morphism or a higher morphism.

%A SIMPLEX IS NOT (A MODEL OF) A HIGHER MORPHISM IN ANY SENSE

\begin{warn}
    In this paper, we still use the term ``categories'' to mean ``$1$-categories'', but sometimes we might use ``ordinary categories'' to emphasize that they are not actually $\infty$-categories. However, we do view $1$-categories as a special kind of $\infty$-categories, namely those with no nontrivial higher morphisms. If one uses quasi-categories to model $\infty$-categories, then ordinary categories are viewed as $\infty$-categories using the nerve functor $\rN(-)$ from \cite[\texttt{002N}, \texttt{003D}]{kerodon}. $1$-morphisms in $\infty$-categories will be simply called \emph{morphisms}, and $k$-morphisms for $k>1$ will be called \emph{higher morphisms} or \emph{higher homotopies}, since they are all invertible. 
\end{warn}

\begin{warn}
    A morphism $f: X \to Y$ between objects $X$ and $Y$ in an $\infty$-category $\cC$ is an isomorphism if it induces an isomorphism in the homotopy $1$-category $\mathrm{h}\cC$ (see \cite[\texttt{004R}]{kerodon}). In such a case, we write $X = Y$ or $X\xra{\sim}Y$ instead of $X \cong Y$. Unlike usual equality, this notation implicitly retains the information of the isomorphism $f$.
        
        % As a result, the meaning of ``isomorphism class" in an $\infty$-category is stronger than that in the $1$-category \cite[\texttt{01D0}, \texttt{01D7}]{kerodon}).
        
    Depending on the $\infty$-category context, isomorphisms can take the form of chain isomorphisms, chain homotopy equivalences, quasi-isomorphisms, or equivalences of $\infty$-categories. We also use $\cC\simeq \cC'$ or $\cC\xra{\sim}\cC'$ to denote the equivalence of two $\infty$-categories. 
    
    Additionally, we use the equality sign to denote natural isomorphisms between two functors, viewed as objects in the $\infty$-category of functors. 
    
    This notation of equality also applies to morphisms and higher morphisms. For example, when we write $h=g\circ f$, we also implicitly record the $2$-morphism that ``witnesses'' the equation (see \cite[\texttt{0042}]{kerodon}). Also, commutative diagrams are witnessed by higher morphisms (see \cite[\texttt{005V}]{kerodon}). %Originally "2-morphisms" here, but actually larger diagrams might be witnessed by higher morphisms, for example for $\bbox^n$-diagrams, there will be n-morphisms witnessing them. 
\end{warn}

To state our main propositions, we introduce the following definitions. We choose to keep the setup as simple as possible while retaining the main ideas; many results in \S \ref{sec: algebraic results} will in fact hold in greater generality.

\bdefn\label{defn: intro}
    \benum
        \item Given $n\in\N$, let $\bbox^n=(\Delta^1)^n$ be the \emph{simplicial cube of dimension $n$}, where $\Delta^1=\{0\to 1\}$ is the category with two objects $0$ and $1$, two identity morphisms, and one morphism from $0$ to $1$. For $n=0$, $\bbox^0$ is just the category with one object and no morphisms except the identity. Note that $\Fun(\bbox^0,\cC)=\cC$ for any $\infty$-category $\cC$.
        % \item For $r\in\N$, let $\BN^r$ be the classifying category of the monoid $\N^r$, so it has only one object whose endomorphism monoid is $\N^r$. Note that $\BN=  \Delta^1/\partial \Delta^1$ from \cite[\texttt{00JC}]{kerodon}, so $\BN^r=  (\Delta^1/\partial \Delta^1)^r$.
        \item Let $\bK$ be a commutative Noetherian ring of finite global dimension in the sense of Definition \ref{defn: global dimension} (see also \cite[\texttt{00O6}]{stacks}), and let $m$ be a positive integer. Assume that either $m$ is even or $\bK$ has characteristic $2$. Examples of such rings $\bK$ include $\Z$, fields, polynomial rings and formal power series rings in finitely many variables over such rings, and localizations of such rings. Following \cite[Definition 3.51]{Ye2026infty}, let \[\De_m(\bK)=\Ndg_\bullet(\Ch_m(\Mod(\bK)_{\proj}))\] be the derived $\infty$-category of (possibly infinitely presented) $m$-periodic chain complexes over $\bK$ and let \[\De_m^c(\bK)=\Ndg_\bullet(\Ch_m(\Mod(\bK)_{\proj}^\fin))\] be the full sub-$\infty$-category of $\De_m(\bK)$ consisting of complexes of finitely presented projective modules, which are exactly compact objects. In the notation of \cite[Example 3.58]{Ye2026infty}, $\De_m^c(\bK)$ is written as $\De_m(\Mod(\bK)^\fin)$. We write $\De(\bK)$ for the usual (non-periodic and unbounded) derived $\infty$-category constructed in \cite[Definition 1.3.5.8]{HAlurie2017} and define $\De^c(\bK)$ to be the full sub-$\infty$-category of compact objects, which are bounded complexes of finitely presented projective modules.
        %We do not use $\De_0(\bK)$ to denote $\De(\bK)$ because some properties of $\De(\bK)$ and $\De_m(\bK)$ are different. 
        From \cite[Proposition 1.3.2.10 and Corollary 1.3.2.18]{HAlurie2017} and \cite[\S 3.6]{Ye2026infty}, all those derived categories are stable in the sense of \cite[Definition 1.1.1.9]{HAlurie2017} and also symmetric monoidal. The compact sub-$\infty$-categories $\De_m^c(\bK)$ and $\De^c(\bK)$ are furthermore rigid. One of the main advantages of stable $\infty$-categories over triangulated categories is that the mapping cone is the cofiber, which is a particular kind of colimit and is hence functorial. 
        %Moreover, if either $m$ is even or $2=0$ in $\bK$, then all those derived categories are symmetric monoidal with respect to the derived tensor product. (NOT CLEAR. ON THE OTHER HAND, I THINK THAT WE DO NOT USE TENSOR PRODUCTS OF PERIODIC COMPLEXES.) 
        \item Given $r\in\N$, let $A_r=\bK[T_1,\ldots,T_r]$ be the polynomial algebra over $\bK$. Given $n\in\N$ and $m\in\Z_+$, let $\De$ be one of the following $\infty$-categories
        % \begin{equation}\label{eq: derived cat}
        %         \De^{n}(\bK)=\Fun(\bbox^n,\De(\bK)),~\De_{m}^{n}(\bK)=\Fun(\bbox^n,\De_m(\bK)),
        %         \end{equation}
                \begin{equation}\label{eq: derived cat compact}
            \De^{n,c}(\bK)=\Fun(\bbox^n,\De^c(\bK)),\aand\De_{m}^{n,c}(\bK)=\Fun(\bbox^n,\De_m^c(\bK)).
        \end{equation}
        Let $\End^r(\De)=\Fun(\BN^r,\De)$, where $\BN^r$ is the classifying $\infty$-category of the monoid $\N^r$ from \cite[\texttt{04G0}]{kerodon}. We also write $\Mod_{A_r}(\De)$ for $\End^r(\De)$, viewed as the full subcategory of $A_r$-modules in the corresponding ambient derived category whose underlying objects lie in $\De$. An object in $\End^r(\De)$ is written as $\ov{\sX}=(\sX,T)$ (also $(X,T)$ in \S \ref{subsec: Eigenspace decomposition}), where $\sX\in \De$ and $T$ consists of $r$ endomorphisms of $\sX$ that commute homotopy-coherently.
        \item Given $A_r$ and $\De$ as above and an ideal $I$ of $A_r$, let $\End^r_{I}(\De)=\Mod_{A_r}^{\Nil(I)}(\De)$ be the full sub-$\infty$-category of $\End^r(\De)$ consisting of $I$-nilpotent objects as in \cite[Remark 7.0.0.1]{SAGlurie2018}, i.e.\ objects $(\sX,T)$ such that for each $f\in I$, 
        \begin{equation}\label{eq: colimit zero}
            \colim(\sX\xrightarrow{f}\sX\xrightarrow{f}\cdots)=0,
        \end{equation}
        where $\mathrm{colim}$ is the homotopy colimit in $\Fun(\bbox^n,\De(\bK))$ or $\Fun(\bbox^n,\De_m(\bK))$, respectively. In the ambient $A_r$-module $\infty$-category, where the underlying objects need not be compact, let $\RG_I$ denote the \emph{local cohomology functor}, the right adjoint to the inclusion of the full subcategory of $I$-nilpotent objects (see \cite[Definition 7.1.2.1]{SAGlurie2018} and \cite[\texttt{0952}]{stacks}).
        \item Suppose that $\bK$ is algebraically closed. By Hilbert's Nullstellensatz, every maximal ideal of $A_r=\bK[T_1,\ldots,T_r]$ has the form
\begin{equation}\label{eq: maximal ideal}
    \fm_{\vec{\la}}=(T_1-\la_1,\ldots,T_r-\la_r)
\end{equation}
for a unique $\vec{\la}=(\la_1,\ldots,\la_r)\in\bK^r$. For $\ov{\sX}\in\End^r(\De)$, write $\ov{\sX}_{\vec{\la}}=\RG_{\fm_{\vec{\la}}}(\ov{\sX})$ and denote its underlying object by $\sX_{\vec{\la}}$.
        % \item Given $A$ and $\De$ as above and a finitely presented ideal $I=(f_1,\ldots,f_k)$ of $A$, for any object $(X,T)$ of $\De$, define 
        % \begin{equation}\label{eq: total fiber}
        %     \RG_I(X,T)=\colim_{N\to\infty}\Tfib(f_1^N,\ldots,f_k^N: X\to X),
        % \end{equation}
        % where $\Tfib$ denotes the total fiber, i.e.\ taking the fiber of $f_1^N: A\to A$, then taking the fiber of $f_2^N$ acting on the previous fiber, and so on, taking fibers $k$ times. See also Construction \ref{cons: total cofiber} for the discussion on the total cofiber. Note that when $\bK=\C$ or other algebraically closed field, then by Hilbert's Nullstellensatz, any maximal ideal $\fm$ of $A=\bK[T_1,\ldots,T_r]$ has the form $\fm=(T_1-\la_1,\ldots,T_r-\la_r)$ for $\la_1,\ldots,\la_r\in \bK$. When $I=\fm$ is a maximal ideal of $A$, the colimit of the total fibers (\ref{eq: total fiber}) can be thought of as the ``union'' of elements killed by higher and higher powers of $T_i-\la_i$, and is the analog of the simultaneous generalized eigenspace with eigenvalue $(\la_1,\ldots,\la_r)$. Note that from the formula (\ref{eq: total fiber}), $\RG_I$ seems to depend on the choice of generators $f_1,\ldots,f_k$ of $I$, but we will see in \S \ref{subsec: Eigenspace decomposition} that it actually does not, by rewriting it categorically as a left adjoint. 
    \eenu
\edefn

The following theorem provides a dictionary between the original definition of hypercubes and the definition using $\infty$-categories, which essentially follows from \cite[Constructions 7.6 and 7.7]{Ye2026infty} and Construction \ref{concrete-cell-unstable}.

\bthm\label{thm: main infty replacement}
   Let $\De$ be as in \eqref{eq: derived cat compact}. The objects of $\De$ and $\End^r(\De)$ admit the following descriptions.
    \benum
        \item An object in $\De$ corresponds to a hypercube of dimension $n$ in (\ref{eq: hypercube}) with either $\Z$ or $\Z/m$-homological grading;
        \item An object in $\End^r(\De)$ corresponds to a hypercube of dimension $n$ together with either $\Z$ or $\Z/m$-homological grading and $r$ homotopy-coherently commutative operators on the hypercube.
    \eenu
\ethm
\begin{notation}
    In \cite{Ye2026infty}, the notion in the $\infty$-categorical setting is called a \emph{cube complex} for distinction. In view of Theorem \ref{thm: main infty replacement}, we use the terms hypercube and cube complex interchangeably.
\end{notation}

Then we provide the eigendecomposition of hypercubes via $\infty$-categories. The proof of the following theorem will be given in \S \ref{subsec: Eigenspace decomposition}. 

\bthm\label{thm: main eigenspace decomposition}
Let $\bK$ be an algebraically closed field and let $\De$ be either $\De^c(\bK)$ or $\De_m^c(\bK)$. The direct sum functor gives an equivalence of stable $\infty$-categories
\begin{equation}\label{eq: main eigendecomposition}
\bigoplus_{\vec{\la}\in\bK^r}\End^r_{\fm_{\vec{\la}}}(\De)\xra{\sim}\End^r(\De),
\end{equation}
where the direct sum consists of families with finite support. Its inverse sends $\ov{\sX}$ to $(\RG_{\fm_{\vec{\la}}}(\ov{\sX}))_{\vec{\la}\in\bK^r}$. With the symmetric monoidal structures from Remarks \ref{rem: periodic-complexes-as-dg-modules} and \ref{rem: periodic-monoidal-decomposition}, this equivalence is symmetric monoidal.
\ethm

\bcor\label{cor: main eigen}
Let $\De$ be as in Theorem \ref{thm: main eigenspace decomposition}. There is an equivalence of stable $\infty$-categories
\begin{equation}\label{eq: eigen for hypercube}
\End^r(\Fun(\bbox^n,\De))\simeq\bigoplus_{\vec{\la}\in\bK^r}\Fun(\bbox^n,\End^r_{\fm_{\vec{\la}}}(\De)).
\end{equation}
This equivalence is symmetric monoidal for the structures of Theorem \ref{thm: main eigenspace decomposition} applied pointwise in the cube direction.
\ecor
\bpf
By Theorem \ref{thm: main eigenspace decomposition}, we have
\[
\begin{aligned}
\End^r(\Fun(\bbox^n,\De))
&=\Fun(\BN^r\times\bbox^n,\De)\\
&=\Fun(\bbox^n,\End^r(\De))\\
&\simeq\Fun\Bigl(\bbox^n,\bigoplus_{\vec{\la}\in\bK^r}\End^r_{\fm_{\vec{\la}}}(\De)\Bigr)\\
&=\bigoplus_{\vec{\la}\in\bK^r}\Fun(\bbox^n,\End^r_{\fm_{\vec{\la}}}(\De)).
\end{aligned}
\]
The last identification holds because $\bbox^n$ has finitely many vertices, so the union of their supports is finite. Symmetric monoidality follows by applying Theorem \ref{thm: main eigenspace decomposition} pointwise in the cube direction.
\epf

\bdefn\label{defn: eigen functor}
The equivalence \eqref{eq: eigen for hypercube}, followed by forgetting the operators on each summand, gives the \emph{eigendecomposition functor}
\[\Eig^{n,r}:\End^r(\Fun(\bbox^n,\De))\to\bigoplus_{\vec{\la}\in\bK^r}\Fun(\bbox^n,\De).\]
We write $\Eig^{n,r}_{\vec{\la}}(\ov{\sX})=\sX_{\vec{\la}}$ for its component at $\vec{\la}$, called the \emph{generalized eigenspace associated to $\vec{\la}$}. Note that on homology they give the usual generalized eigenspaces.
\edefn

Without passing to the $\infty$-categorical framework, one needs to choose a splitting map $\Phi$ in \eqref{eq: map phi} to obtain the hypercube of the eigendecomposition. By \cite[Lemma 6.5]{GZ2021}, the splitting map can be chosen inductively over subcubes so that the resulting hypercube behaves well with respect to subcubes. However, the proof of that lemma still involves some arbitrary choices. The results of \cite[Lemma 6.7]{GZ2021} and its generalization to hypercubes imply that different choices of the splitting map induce homotopy equivalent hypercubes.

Moreover, for two homotopy-coherently commutative operators $T_1$ and $T_2$, different orders of applying the eigendecomposition give homotopy equivalent hypercubes by \cite[Lemma 6.8]{GZ2021}. By \cite[Proposition 6.9]{GZ2021}, if $T_2-T_1=\al\id$ for some $\al\in\C$, then the corresponding eigendecompositions differ by a shift of eigenvalues up to homotopy equivalence. The following form of this comparison allows $T_2-T_1-\al\id$ to be nilpotent and is used in \S \ref{sec: mixed hypercube}.

\bcor\label{cor: nilpotent}
Let $\De$ be as in Theorem \ref{thm: main eigenspace decomposition} and let $\ov{\sX}=(\sX,T)\in\End^r(\Fun(\bbox^n,\De))$. Fix $i,j\in[1,r]$ and $\al\in\bK$. If $T_j-T_i-\al\id$ acts nilpotently on the homology of every vertex, then $\sX_{\vec{\la}}=0$ unless $\lambda_j=\lambda_i+\al$. In particular, the generalized eigenspaces of $T_i$ at $\lambda$ and of $T_j$ at $\lambda+\al$ are naturally identified by the inclusions of their joint generalized eigenspace.
\ecor
\bpf
On the homology of each vertex of $\sX_{\vec{\la}}$, the three terms on the right of
\[(\lambda_j-\lambda_i-\al)\id=(T_j-T_i-\al\id)-(T_j-\lambda_j\id)+(T_i-\lambda_i\id)\]
commute and are nilpotent. If $\lambda_j-\lambda_i-\al\ne0$, the homology vanishes at every vertex, so $\sX_{\vec{\la}}=0$. The last assertion follows by summing the joint blocks with the prescribed eigenvalues.
\epf
\subsection{Combination of bypass maps}\label{subsec: combination of bypass maps}
In this subsection, we introduce the setup for combining bypass maps in future work by the second author and collaborators. We start with the simplest case corresponding to $n=k=1$ in Theorem \ref{thm: main bypass cube} and then consider the general case.

Suppose $(M,\ga_0)$ is a balanced sutured manifold and $\eta\subset \partial M$ is an oriented simple closed curve intersecting the suture $\ga_0$ twice. There are two natural bypass arcs associated to $\eta$, as shown in Figure \ref{fig: two bypass arcs}, called the \emph{positive} and \emph{negative} bypass arcs $\eta^+$ and $\eta^-$, respectively. The induced sutures $\ga_1^\pm$ in Figure \ref{fig:bypass} for these two bypass arcs are isotopic, so we omit the sign in the superscript. The top sutures $\ga_\pm$ in Figure \ref{fig:bypass} for $\eta^\pm$ are generally not isotopic, with some exceptions, such as when $\partial M\cong T^2$ and $\ga$ consists of two parallel simple closed curves.
\definecolor{zlgreen}{RGB}{1,102,70}

\newcommand{\ZLLeftArrow}[4]{%
  \draw[#1,-{Stealth[length=2.2mm,width=1.5mm]}] (#2,#4)--(#3,#4);%
}

\newcommand{\ZLTwoSimpleTube}[2]{%
  \begin{scope}[shift={(#1,#2)}]
    \draw (0,0) ellipse[x radius=\zlTwoR,y radius=\zlTwoH];
    \draw (0,\zlTwoH)--(\zlTwoL,\zlTwoH) (0,-\zlTwoH)--(\zlTwoL,-\zlTwoH);
    \draw[dashed] (\zlTwoL,\zlTwoH) arc[start angle=90,end angle=270,x radius=\zlTwoR,y radius=\zlTwoH];
    \draw (\zlTwoL,-\zlTwoH) arc[start angle=-90,end angle=90,x radius=\zlTwoR,y radius=\zlTwoH];
  \end{scope}%
}

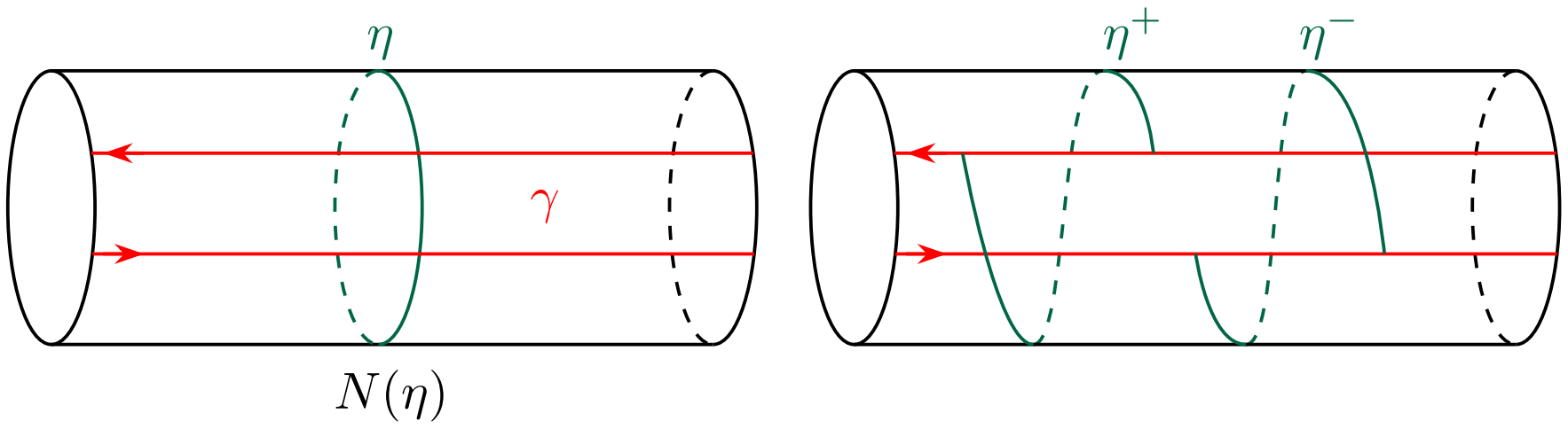
\begin{figure}[htbp]
\centering
\pgfmathsetlengthmacro{\zlUnit}{0.8\textwidth/24}
\begin{tikzpicture}[
  x=\zlUnit,y=\zlUnit,line width=0.25mm,line cap=butt,line join=round,
  every node/.style={font=\normalsize,anchor=base west,inner sep=0pt}
]
\def\zlTwoL{9.71}
\def\zlTwoR{0.64}
\def\zlTwoH{2.01}
\path[use as bounding box] (0,0) rectangle (24,6);

\ZLTwoSimpleTube{1.40}{2.95}
\begin{scope}[shift={(6.20,2.95)}]
  \draw[zlgreen,dashed] (0,\zlTwoH) arc[start angle=90,end angle=270,x radius=\zlTwoR,y radius=\zlTwoH];
  \draw[zlgreen] (0,-\zlTwoH) arc[start angle=-90,end angle=90,x radius=\zlTwoR,y radius=\zlTwoH];
\end{scope}

\begin{scope}[shift={(13.18,2.95)}]
  \draw[zlgreen,dashed]
    (2.63,-\zlTwoH)..controls (3.14,-2.00) and (3.00,2.00)..(3.65,\zlTwoH)
    (5.75,-\zlTwoH)..controls (6.27,-2.00) and (6.09,2.00)..(6.63,\zlTwoH);
\end{scope}
\ZLTwoSimpleTube{13.18}{2.95}

\foreach \xx in {1.40,13.18}{%
  \begin{scope}[shift={(\xx,2.95)}]
    \foreach \yy in {0.80,-0.68}{%
      \draw[red] ({\zlTwoR*sqrt(1-(\yy/\zlTwoH)^2)},\yy) --
        ({\zlTwoL+\zlTwoR*sqrt(1-(\yy/\zlTwoH)^2)},\yy);
    }%
    \ZLLeftArrow{red}{1.36}{0.76}{0.80}
    \ZLLeftArrow{red}{0.76}{1.36}{-0.68}
  \end{scope}%
}

\begin{scope}[shift={(13.18,2.95)}]
  \draw[zlgreen]
    (1.59,0.80)..controls (1.59,0.80) and (2.09,-2.02)..(2.63,-\zlTwoH)
    (3.65,\zlTwoH)..controls (4.29,2.02) and (4.39,0.80)..(4.39,0.80)
    (5.01,-0.68)..controls (5.01,-0.68) and (5.19,-2.01)..(5.75,-\zlTwoH)
    (6.63,\zlTwoH)..controls (7.53,2.01) and (7.78,-0.68)..(7.78,-0.68);
\end{scope}

\node at (5.52,0) {$N(\eta)$};
\node[zlgreen] at (6.00,5.28) {$\eta$};
\node[red] at (8.40,2.88) {$\gamma$};
\node[zlgreen] at (16.80,5.28) {$\eta^+$};
\node[zlgreen] at (19.68,5.28) {$\eta^-$};
\end{tikzpicture}
\caption{The curve $\eta$ and its two bypass arcs $\eta^+$ and $\eta^-$. On the visible side, the outward normal to $\partial M$ points into the cylinder.}
\label{fig: two bypass arcs}
\end{figure}

Let $S\subset (M,\ga_0)$ be an admissible surface, possibly intersecting $\eta^+$ and $\eta^-$, such that the algebraic intersection number \[a=\partial S\cdot \eta\]is positive. Let $S^+$ and $S^-$ be obtained from $S$ by isotoping its boundary in $N(\eta)$ as in Figure \ref{fig:isotopy} so that\[S^+\cap \eta^+=\emptyset\aand S^-\cap \eta^-=\emptyset,\]and the induced surfaces in $(M,\ga_1)$ are isotopic when fixing $\ga_1$. Note that the operation of isotoping the boundary is referred to as \emph{stabilization} in Figure \ref{fig:stabilization} and corresponds to a grading shift for the surface grading in Proposition \ref{prop: grading properties SHI} \eqref{g3}. It is straightforward to verify from Figures \ref{fig:isotopy} and \ref{fig:stabilization} that the grading shift between the surface gradings associated to $S^+$ and $S^-$ is exactly $2a$, where the factor of $2$ comes from our rescaling of the grading.

\newcommand{\ZLEightSimpleTube}[2]{%
  \begin{scope}[shift={(#1,#2)}]
    \draw (0,0) ellipse[x radius=\zlEightR,y radius=\zlEightH];
    \draw (0,\zlEightH)--(\zlEightL,\zlEightH) (0,-\zlEightH)--(\zlEightL,-\zlEightH);
    \draw[dashed] (\zlEightL,\zlEightH) arc[start angle=90,end angle=270,x radius=\zlEightR,y radius=\zlEightH];
    \draw (\zlEightL,-\zlEightH) arc[start angle=-90,end angle=90,x radius=\zlEightR,y radius=\zlEightH];
    \foreach \yy in {0.61,-0.98}{%
      \draw[red] ({\zlEightR*sqrt(1-(\yy/\zlEightH)^2)},\yy) --
        ({\zlEightL+\zlEightR*sqrt(1-(\yy/\zlEightH)^2)},\yy);
    }%
    \ZLLeftArrow{red}{1.36}{0.76}{0.61}
    \ZLLeftArrow{red}{0.76}{1.36}{-0.98}
  \end{scope}%
}

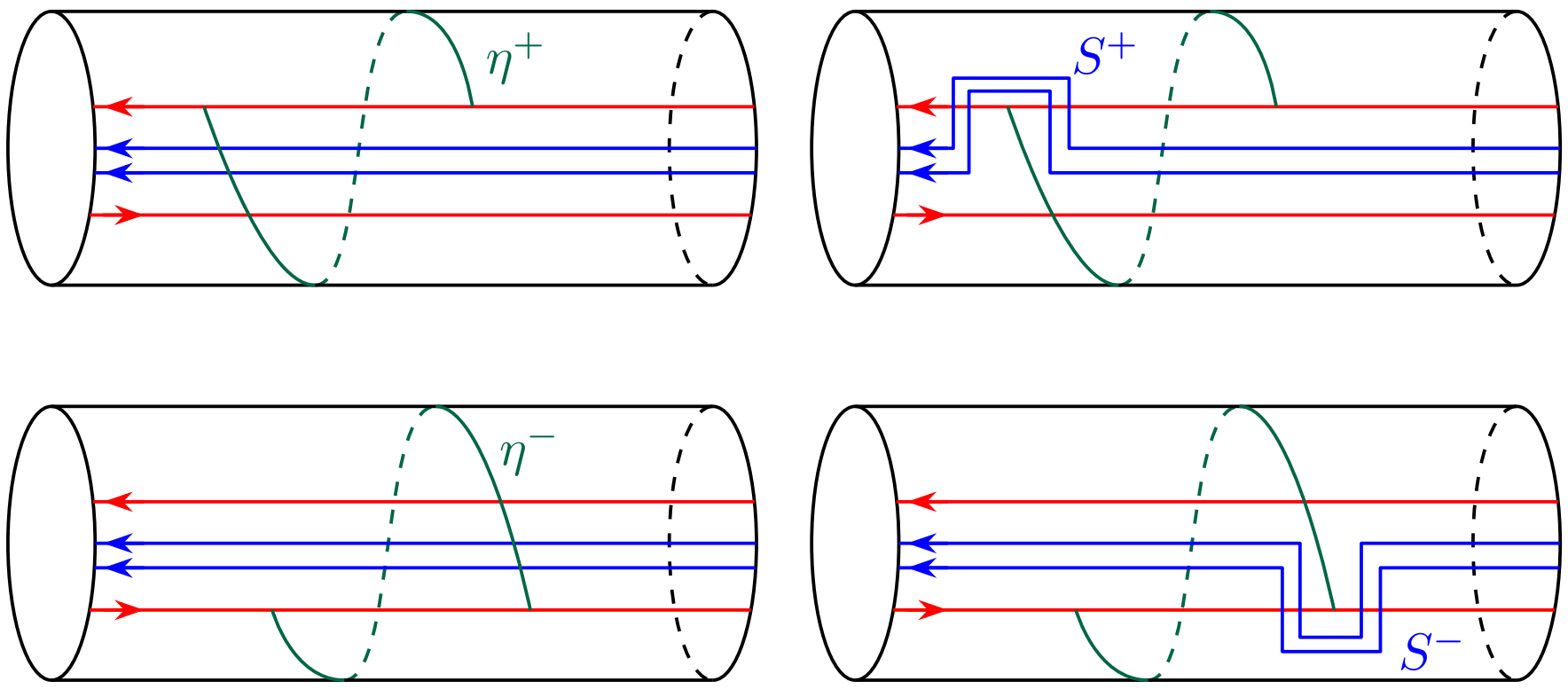
\begin{figure}[htbp]
\centering
\pgfmathsetlengthmacro{\zlUnit}{0.8\textwidth/24}
\begin{tikzpicture}[
  x=\zlUnit,y=\zlUnit,line width=0.25mm,line cap=butt,line join=round,
  every node/.style={font=\normalsize,anchor=base west,inner sep=0pt}
]
\def\zlEightL{9.71}
\def\zlEightR{0.64}
\def\zlEightH{2.01}
\path[use as bounding box] (0,0) rectangle (24,12.75);

\foreach \xx in {1.35,13.15}{%
  \begin{scope}[shift={(\xx,9.70)}]
    \draw[zlgreen,dashed] (3.87,-\zlEightH)..controls (4.63,-2.00) and (4.38,2.00)..(5.21,\zlEightH);
  \end{scope}
  \begin{scope}[shift={(\xx,3.90)}]
    \draw[zlgreen,dashed] (4.30,-\zlEightH)..controls (5.06,-2.00) and (4.81,2.00)..(5.64,\zlEightH);
  \end{scope}
}

\foreach \xx in {1.35,13.15}{%
  \foreach \yy in {3.90,9.70}{\ZLEightSimpleTube{\xx}{\yy}}%
}

\foreach \yy in {3.90,9.70}{%
  \begin{scope}[shift={(1.35,\yy)}]
    \foreach \railY in {0,-0.36}{%
      \draw[blue] ({\zlEightR*sqrt(1-(\railY/\zlEightH)^2)},\railY) --
        ({\zlEightL+\zlEightR*sqrt(1-(\railY/\zlEightH)^2)},\railY);
    }%
  \end{scope}%
}
\begin{scope}[shift={(13.15,9.70)}]
  \draw[blue]
    (0.64,0)--(1.44,0)--(1.44,1.03)--(3.14,1.03)--(3.14,0)--(10.35,0)
    (0.63,-0.36)--(1.67,-0.36)--(1.67,0.84)--(2.86,0.84)--(2.86,-0.36)--(10.34,-0.36);
\end{scope}
\begin{scope}[shift={(13.15,3.90)}]
  \draw[blue]
    (0.64,0)--(6.53,0)--(6.53,-1.38)--(7.43,-1.38)--(7.43,0)--(10.35,0)
    (0.63,-0.36)--(6.27,-0.36)--(6.27,-1.59)--(7.71,-1.59)--(7.71,-0.36)--(10.34,-0.36);
\end{scope}

\foreach \xx in {1.35,13.15}{%
  \foreach \yy in {3.90,9.70}{%
    \begin{scope}[shift={(\xx,\yy)}]
      \ZLLeftArrow{blue}{1.36}{0.76}{0}
      \ZLLeftArrow{blue}{1.36}{0.76}{-0.36}
    \end{scope}%
  }%
}

\foreach \xx in {1.35,13.15}{%
  \begin{scope}[shift={(\xx,9.70)}]
    \draw[zlgreen]
      (2.24,0.61)..controls (2.24,0.61) and (3.09,-2.02)..(3.87,-\zlEightH)
      (5.21,\zlEightH)..controls (5.99,2.02) and (6.18,0.61)..(6.18,0.61);
  \end{scope}
  \begin{scope}[shift={(\xx,3.90)}]
    \draw[zlgreen]
      (3.24,-0.98)..controls (3.24,-0.98) and (3.53,-2.01)..(4.30,-\zlEightH)
      (5.64,\zlEightH)..controls (6.42,2.02) and (7.03,-0.98)..(7.03,-0.98);
  \end{scope}
}

\node[zlgreen] at (7.9,5.05) {$\eta^-$};
\node[zlgreen] at (7.7,10.80) {$\eta^+$};
\node[blue] at (16.3,10.80) {$S^+$};
\node[blue] at (21.1,2.05) {$S^-$};
\end{tikzpicture}
\caption{The oriented boundary of $S$ is shown in blue on the left, with $\eta^+$ above and $\eta^-$ below. On the right, $S^+$ is obtained by positive stabilizations and $S^-$ by negative stabilizations. We use the orientation convention of Figure \ref{fig: two bypass arcs}.}
\label{fig:isotopy}
\end{figure}
We choose $S^\pm$ to be admissible, applying the same stabilization to both away from $N(\eta)$ if necessary. We then apply Theorem \ref{thm: main bypass cube} to $(S^+,\eta^+)$ and $(S^-,\eta^-)$ separately. Since the spectral sequence in Theorem \ref{thm: main bypass cube} is from the cube filtration, similar to \eqref{eq: cone}, we indeed obtain the quasi-isomorphisms\begin{equation}\label{eq: bypass cone}
    SHI(M,\ga_\pm,S^\pm,\lambda)\cong \cone\big(SHI(M,\ga_0,S^\pm,\lambda)\to SHI(M,\ga_1,S^\pm,\lambda)\big),
\end{equation}where we use the fact that a chain complex over $\C$ is quasi-isomorphic to its homology, and the notation $S^\pm$ and $\lambda\in 2\Z$ is for the surface grading associated. We denote the maps in the mapping cones of \eqref{eq: bypass cone} by $\psi_\pm$ and call them the \emph{positive} and \emph{negative} bypass maps. The previous assumption on $S^\pm$ implies that \begin{equation}\label{eq: shift grading}
\begin{aligned}
        SHI(M,\ga_0,S^+,\lambda)&\cong SHI(M,\ga_0,S^-,\lambda-2a),\\\aand SHI(M,\ga_1,S^+,\lambda)&\cong SHI(M,\ga_1,S^-,\lambda).
\end{aligned}
\end{equation}Hence we can construct a mapping cone\begin{equation}\label{eq: mixed mapping cone}
    \cone\big(\bigoplus_{\lambda\in 2\Z}SHI(M,\ga_0,S^+,\lambda)\xra{\psi_++\psi_-}\bigoplus_{\lambda\in 2\Z}SHI(M,\ga_1,S^+,\lambda)\big),
\end{equation}called the \emph{mixed mapping cone}, such that $\psi_+$ preserves the grading $\lambda$ and $\psi_-$ decreases the grading by $2a$. Then the grading induces a filtration on \eqref{eq: mixed mapping cone} such that the $E_1$-page is isomorphic to \[\bigoplus_{\lambda\in 2\Z}SHI(M,\ga_+,S^+,\lambda)\] via \eqref{eq: bypass cone}. We refer to such a filtration as a \emph{grading filtration}, in contrast to the \emph{cube filtration} in the introduction. Although we can also take another filtration from the grading so that the $E_1$-page is isomorphic to $\bigoplus_{\lambda\in 2\Z}SHI(M,\ga_-,S^-,\lambda)$, we do not need this one in this paper.

The construction of the mixed mapping cone generalizes easily for multiple admissible surfaces $S_l\subset (M,\ga_0)$ for $l=1,\dots,k$. With respect to the surface multi-grading associated to $S_l^+$, the map $\psi_+$ has degree zero and $\psi_-$ has multi-degree $-2\vec{a}$, where $\vec{a}=(a_1,\dots,a_k)$ and $a_l=\pa S_l\cdot\eta$.

The generalization to several bypass arcs is more complicated. Suppose now that $\eta_i\subset \partial M$, $i=1,\dots,n$, are pairwise disjoint oriented simple closed curves, each intersecting $\ga_0$ twice, with pairwise disjoint neighborhoods $N(\eta_i)$. Let $\eta_{i}^\pm$ be the corresponding bypass arcs from Figure \ref{fig: two bypass arcs} and let $\psi_{i,\pm}$ be the corresponding bypass maps. As for $\ga_\pm$ when $n=1$, for $s\in\{+,-\}^n$, let $\ga_s$ be the suture obtained from $\ga_0$ by taking the top suture for $\eta_{i}^{s_i}$ in Figure \ref{fig:bypass}, where $s_i$ is the $i$-th entry of $s$. Also for $\ep\in \{0,1\}^n$, let $\ga_\ep$ denote one of the bottom sutures in Figure \ref{fig:bypass} for $\eta_{i}^\pm$ with any choice of signs, as they are all isotopic. In particular, we have $\ga_{\vec{0}}=\ga_0$.

For any admissible surface $S_l\subset (M,\ga_0)$ with $l=1,\dots,k$, we choose $2^n$ admissible surfaces $S^{s}_l$ for $s\in\{+,-\}^n$ obtained by stabilizations, such that in $N(\eta_i)$, the boundary of $S_l$ is isotopic to the boundary in one of the cases in Figure \ref{fig:isotopy}, depending on the sign $s_i$. For $i=1,\dots,n$, let \[f_i=-2\cdot (\partial S_1\cdot \eta_i,\dots, \partial S_k\cdot \eta_i)\in (2\Z)^k,\]which is the grading shift of $\psi_{i,-}$ when $\psi_{i,+}$ is grading preserving.

Again by Theorem \ref{thm: main bypass cube}, we obtain quasi-isomorphisms for all $s\in \{+,-\}^n$ as in \eqref{eq: iterated cone}\begin{equation}\label{eq: many iterated cones}
    \begin{aligned}
SHI(M,\ga_s,\{S_l^s\}_{l\in[1,k]},\vec{\lambda})&\\[-1ex]
&\hspace{-3em}\cong\left(\bigoplus_{\ep\in\{0,1\}^n}SHI(M,\ga_\ep,\{S_l^s\}_{l\in[1,k]},\vec{\lambda})[|\ep|-n],~\sum_{\ep\leq\epp}D^s_{\ep,\epp}\right).
\end{aligned}
\end{equation}Here we implicitly use the fact that the chain complex over a field is quasi-isomorphic to its homology. See also the explicit construction of the homological perturbation lemma for hypercubes in \cite[Lemma 2.11]{HHSZ22perturb} and \cite[Lemma 3.3]{Josh2023trefoil}.

Generalizing \eqref{eq: shift grading}, we have\begin{equation}\label{eq: many shift grading}
    SHI(M,\ga_{\ep},\{S_l^{\vec{+}}\}_{l\in [1,k]},\vec{\lambda})\cong SHI(M,\ga_{\ep},\{S_l^s\}_{l\in [1,k]},\vec{\lambda}+\sum_{\substack{\ep_i=0\\s_i=-}}f_i),
\end{equation}where $\vec{+}=(+,\dots,+)$.

Moreover, we have $D^s_{\ep,\ep}=0$, and if $|\ep'-\ep|=1$, then the map $D^s_{\ep,\epp}$ depends only on the entry of $s$ corresponding to the coordinate in which $\ep$ and $\epp$ differ, because it is the bypass map on the homology level.

To generalize the mixed mapping cone \eqref{eq: mixed mapping cone} to a mixed hypercube, we need the following assumption.
\begin{assumption}\label{assump: mixed condition}
There exist choices as in Remark \ref{rem: dependence of data} and representatives of the hypercubes \eqref{eq: many iterated cones} in $\Fun(\bbox^n,\De^c_2(\C))$ such that, for any $\ep\le\epp\in\{0,1\}^n$, the map $D^s_{\ep,\epp}$ depends only on the entries of $s$ corresponding to the coordinates in which $\ep$ and $\epp$ differ. We refer to this condition as the \emph{mixed condition}.
\end{assumption}
\brem\label{rem: trivial for n=2}
Assumption \ref{assump: mixed condition} is trivial when $\ep=\vec{0}=(0,\dots,0)$ and $\epp=\vec{1}=(1,\dots,1)$. Since it also holds whenever $|\epp-\ep|\leq 1$, the mixed condition holds for $n = 2$ once the vertex and edge representatives have been chosen in common.
\erem
If Assumption \ref{assump: mixed condition} holds, we can define the sum \begin{equation}\label{eq: mixed differential}
    D^{\mix}=\sum_{\ep\le \epp\in\{0,1\}^n}\sum_{s\in\{+,-\}^{\epp-\ep}}D^s_{\ep,\epp},
\end{equation}where $\{+,-\}^{\epp-\ep}$ consists of those $s\in\{+,-\}^n$ satisfying $s_i=+$ whenever $\ep_i=\ep^\prime_i$. The fact that \[D^s=\sum_{\ep\le \epp\in \{0,1\}^n}D^s_{\ep,\epp}\]is a differential for all $s\in\{+,-\}^n$ then implies that $D^{\mix}$ is also a differential on \[
\bigoplus_{\vec{\lambda}\in(2\Z)^k}\bigoplus_{\ep\in\{0,1\}^n}SHI(M,\ga_\ep,\{S_l^{\vec{+}}\}_{l\in[1,k]},\vec{\lambda})[|\ep|-n].
\]We call the corresponding chain complex a \emph{mixed hypercube}.

Suppose there is a homomorphism $\rho:(2\Z)^k\to\Z$ such that $\rho(f_i)>0$ for all $i\in[1,n]$. The summands with $\rho(\vec{\lambda})\ge p$ define a decreasing filtration on the mixed hypercube. Its $E_1$-page is isomorphic to \[\bigoplus_{\vec{\lambda}\in (2\Z)^k}SHI(M,\ga_{\vec{+}},\{S_l^{\vec{+}}\}_{l\in [1,k]},\vec{\lambda})\]via \eqref{eq: many iterated cones} for $s=\vec{+}$. We call such a filtration a \emph{$\rho$-grading filtration}. The mixed hypercube and the $\rho$-grading filtration are the desired constructions mentioned in Remark \ref{rem: further direction}.

Now we focus on Assumption \ref{assump: mixed condition}, which is the main motivation of this paper. We need to unpack the definitions of $SHI$ and the surface multi-grading. Namely, we consider the closures of balanced sutured manifolds and the extensions of the admissible surfaces to the closures, as described in the introduction. Since we do not discuss the details of the topological construction of closures in this paper, we summarize the construction as a topological assumption in the statement of the theorem and will verify it in future work by the second author and collaborators.

As a heuristic, we start with the topological assumption for $n=1$ and general $k$. From the sketch in the introduction and the details in \S \ref{sec: Link spectral sequences for sutured instanton homology}, to define the surface multi-grading, we need to find a \emph{special} closure\footnote{Indeed, the closure is an abbreviation for the marked odd closure in Definition \ref{defn: closure}, which also includes the information of the bundle set $\omega$. Since the bundle set behaves similarly to $R$ in further constructions, such as surgery cobordisms and canonical cobordisms, it is safe to omit it from the notation here and consider it carefully in \S \ref{sec: mixed hypercube}.} $(Y_0,R)$ of $(M,\ga_0)$, such that the following holds.
\begin{enumerate}
    \item The admissible surfaces $S_l^\pm$ for $l\in[1,k]$ extend to closed surfaces $\ov{S}_l^\pm$ in $Y_0$, necessarily intersecting $R$. We do not need the extensions of $S_l$.
    \item The positive and negative bypass maps are induced by surgery cobordisms $W^\pm$ with source $Y_0$ along curves $C^\pm$ disjoint from $R$. Moreover, we assume that $C^+$ is disjoint from $\ov{S}_l^+$, and $C^-$ is disjoint from $\ov{S}_l^-$. Let $Y^\pm_1$ be the target of $W^\pm$ and let $R$ also denote the resulting surface in $Y^\pm_1$. Also, let $\ov{S}_l^\pm\subset Y^\pm_1$ denote the resulting surfaces. Then by construction, both $(Y^\pm_1,R)$ are closures of $(M,\ga_1)$ and $\ov{S}_l^\pm\subset Y^\pm_1 $ induce the multi-grading for $(M,\ga_1)$.
    \item Let $X$ be the canonical cobordism from $(Y^-_1,R)$ to $(Y^+_1,R)$, which is indeed a cobordism \begin{equation}\label{eq: canonical cobordism}
        (X,\wti{Z}_R):(Y^-_1,R)\sqcup (Y_\varphi,R_{\varphi})\to (Y^+_1,R\sqcup R_{\varphi})
    \end{equation}induced by Floer's excision along $R$ and $R_{\varphi}$ (cf.\ \cite[\S 3.2]{kronheimer2010knots}), where $Y_\varphi=(R_{\varphi}\times I)/(x,1)\sim(\varphi(x),0)$ is the mapping torus of $\varphi$, and the cobordism\footnote{We use $\wti{Z}_R$ instead of $Z_R$ so that the notation is compatible with \cite[Construction 8.9]{Ye2026infty}, where $Z_R$ is reserved for the surfaces on the boundary of $\wti{Z}_R$, namely copies of $R$.} $\wti{Z}_R$ for surfaces is the disjoint union of the product cobordisms of $R$ and $R_\varphi$. From the construction in Figure \ref{fig:isotopy}, we know that the resulting surfaces of $S_l$ in $(M,\ga_1)$ are isotopic when fixing $\ga_1$, and we need to further relate its extension in different closures $(Y^\pm_1,R)$. The special closure should be chosen so that in the canonical cobordism, there are further cobordisms of surfaces \begin{equation*}\label{eq: Z for S}
        \wti{Z}_{l}:(Y^-_1,\ov{S}_l^-)\sqcup (Y_\varphi,F_{l})\to (Y^+_1,\ov{S}_l^+),
    \end{equation*}where $F_l$ is a collection of tori obtained from a collection of simple closed curves on $R$ that is preserved by $\varphi$.
\end{enumerate}

We now state the topological assumption for general $n$. Since there are multiple surfaces $\ov{S}_l^s$ induced by $\ov{S}_l$ that behave differently under bypass maps associated to different $\eta_i$, the assumption on the special closure is considerably more complicated. The reason we do not consider a single surface from the stabilization of $S_l$ that avoids both positive and negative bypass arcs is that it is even more difficult to relate the resulting surfaces under the canonical cobordisms.
\begin{assumption}\label{assump: topological}
There exists a \emph{special} closure $(Y_{\vec{0}},R)$ of $(M,\ga_0)$ satisfying the following conditions, where $i\in[1,n]$, $l\in[1,k]$, $\ep\in\{0,1\}^n$, $s,t\in\{+,-\}^n$, and $\ep_i$ denotes the $i$-th entry of $\ep$ (similarly for $s,t$).
\begin{enumerate}
    \item\label{as1} The admissible surfaces $S_l^s$ extend to closed surfaces $\ov{S}_l^s$ in $Y_{\vec{0}}$.
    \item\label{as2} The bypass maps $\psi_{i,\pm}$ are induced by surgery cobordisms $W^\pm_i$ with source $Y_{\vec{0}}$ along curves $C^\pm_i$ disjoint from $R$. Moreover, the curve $C^{s_i}_i$ is disjoint from $\ov{S}_l^s$. Let $Y^s_\ep$ be obtained from $Y_{\vec{0}}$ by applying surgeries associated to $\psi_{i,s_i}$ for $\ep_i=1$ and let $R$ also denote the resulting surface in $Y^s_\ep$. By construction, we have $Y^s_{\ep}=Y^{s'}_{\ep}$ if $s$ and $s'$ differ only in the coordinates $i$ for which $\ep_i=0$. Also, let $\ov{S}_l^s\subset Y^s_\ep$ denote the resulting surfaces. Then by construction, all $(Y^s_\ep,R)$ are closures of $(M,\ga_\ep)$ and $\ov{S}_l^s\subset Y^s_\ep$ induce the multi-grading for $(M,\ga_\ep)$.
    \item\label{as3} We also write $Y^\pm_i$ for the target of $W^\pm_i$. Let\begin{equation}\label{eq: canonical cobordism 2}
        (X_i,\wti{Z}_R):(Y^-_i,R)\sqcup (Y_{\varphi_i},R_{\varphi_i})\to (Y^+_i,R\sqcup R_{\varphi_i})
    \end{equation}be the canonical cobordism for different closures as in \eqref{eq: canonical cobordism}. This canonical cobordism also induces ones\begin{equation}\label{eq: induced canonical cob}
        (X_{\ep,i}^s,\wti{Z}_R):(Y^s_\ep,R)\sqcup (Y_{\varphi_i},R_{\varphi_i})\to (Y^{s'}_\ep,R\sqcup R_{\varphi_i})
    \end{equation}for $\ep_i=1$, $s_i=-$, and $s'$ is obtained from $s$ by replacing $s_i$ with $s_i'=+$. The special closure should be chosen so that in the canonical cobordism, there are further cobordisms of surfaces \begin{equation}\label{eq: Z for S general}
        \wti{Z}_{i,l}^t:(Y^-_i,\ov{S}_l^t)\sqcup (Y_{\varphi_i},F_{i,l}^t)\to (Y^+_i,\ov{S}_l^{t'}),
    \end{equation}where $t_i=-$ and $t'$ is obtained from $t$ by replacing $t_i$ with $t_i'=+$, and $F_{i,l}^t$ is a collection of tori obtained from a collection of simple closed curves on $R$ that is preserved by $\varphi_i$. If $t_j=s_j$ for every $j$ with $\ep_j=1$, the surface cobordisms in \eqref{eq: Z for S general} induce surface cobordisms in $\smash{X_{\ep,i}^s}$ from \eqref{eq: induced canonical cob}, written as
\begin{equation}\label{eq: Z for S general 2}
    \wti{Z}_{\ep,i,l}^{s,t}:(Y_\ep^s,\ov{S}_l^t)\sqcup (Y_{\varphi_i},F_{i,l}^t)\to (Y_\ep^{s'},\ov{S}_l^{t'}).
\end{equation}
    \item\label{as4} Finally, we assume that some collections of canonical cobordisms and the cobordisms of surfaces are supported in disjoint regions in the sense of \cite[Definitions 7.8 and 8.18]{Ye2026infty}, which is roughly saying that they are product outside some open subsets. More explicitly, consider the cobordisms
    \begin{equation}\label{eq: composition cobordism}
    \begin{aligned}
        W^{+,\prime}_i=&W^+_i\sqcup \bigsqcup_{j=1}^nI\times Y_{\varphi_j}\aand \\W^{-,\prime}_i=&\big(X_i\sqcup \bigsqcup_{j\neq i}I\times Y_{\varphi_j}\big )\circ \big(W_i^-\sqcup \bigsqcup_{j=1}^nI\times Y_{\varphi_j}\big),
    \end{aligned}
    \end{equation}and the cobordisms of surfaces
\begin{equation}\label{eq: composition of cobordism surfaces}
    \wti{Z}^{t,\prime}_{i,l}=\begin{cases}
    I\times \ov{S}^t_l& t_i=+,\\
      \wti{Z}^t_{i,l}\circ (I\times \ov{S}^{t}_l\sqcup  I\times F_{i,l}^t)&t_i=-. 
  \end{cases}
\end{equation}
    Then for any $s\in\{+,-\}^n$, the cobordisms in \begin{equation}\label{eq: all coborodisms}
        \big\{(W^{s_i,\prime}_i,\{\wti{Z}^{t,\prime}_{i,l}\}_{l\in[1,k],t\in \{+,-\}^n, t_i=s_i})\big\}_{i\in[1,n]}
    \end{equation}are supported in disjoint regions, where each $\wti{Z}^{t,\prime}_{i,l}$ is understood to include the product cylinders $I\times F_{j,l}^t$ for $j\ne i$ with $t_j=-$. We also assume that on each cobordism piece the $3$-dimensional cobordisms of surfaces have no quadruple intersections. These conditions allow us to apply the hypercube constructions in \cite[Construction 7.11 and Theorem 8.21]{Ye2026infty}. See \cite[Example 8.22]{Ye2026infty} for the condition on quadruple intersections in the excision model.
\end{enumerate}
\end{assumption}
\brem
For simplicity, we omit the information concerning the bundle sets in Assumption \eqref{assump: topological}. One may impose additional conditions so that the extra bundle sets appearing in Theorem \ref{thm: main grading for cube} can be cancelled using \cite[Remark 7.19]{Ye2026infty}. This will be addressed in future work by the second author and collaborators.
\erem
\bthm\label{thm: mixed hypercube}
Assumption \ref{assump: topological} implies Assumption \ref{assump: mixed condition}. Then we have the construction of the differential $D^\mix$ in \eqref{eq: mixed differential} for the mixed hypercube.
\ethm

\subsection{Organization}
In \S \ref{sec: algebraic results}, we develop algebraic results for $\infty$-categories. We expect readers of this section to be familiar with the basic concepts and terminology of $\infty$-categories. In \S \ref{subsec: Characterization of hypercubes}, we give an alternative characterization of the simplicial cubes $\bbox^n$, which will be used in the proof of Theorem \ref{thm: main infty replacement}. In \S \ref{subsec: Derived category of periodic chain complexes}, we provide several equivalent definitions of the derived $\infty$-category $\De_m(\cA)$ for an arbitrary abelian category $\cA$. In \S \ref{subsec: Eigenspace decomposition}, we prove the eigendecomposition in Theorem \ref{thm: main eigenspace decomposition}, together with several generalizations obtained along the way.

In \S \ref{sec: Link spectral sequences for sutured instanton homology}, we review the construction of sutured instanton homology and prove Theorems \ref{thm: main grading for cube} and \ref{thm: main bypass cube}, based primarily on Theorem \ref{thm: main eigenspace decomposition}.

In \S \ref{sec: mixed hypercube}, we prove Theorem \ref{thm: mixed hypercube} in several steps.

\section{Algebraic results for stable \texorpdfstring{$\infty$}{∞}-categories}\label{sec: algebraic results}
In this section, we focus on algebraic results for stable $\infty$-categories, applicable to general scenarios. The results in \S \ref{subsec: language infty cate} are special cases.

For an object $X$ of a stable $\infty$-category, we write $X[k]=\Sigma^k X$ for $k\in\Z$, where $\Sigma$ denotes suspension, following \cite[Notation 1.1.2.7]{HAlurie2017}. For homologically graded chain complexes, this gives $(X[k])_j=X_{j-k}$. Note that this differs from the convention $(X[k])_j=X_{j+k}$ used in \eqref{eq: hypercube}, which follows \cite[Translation 1.2.8]{Weibel1994}.
\subsection{Characterization of hypercubes}\label{subsec: Characterization of hypercubes}
In this subsection, we provide an alternative description of the simplicial cubes $\bbox^n=(\Delta^1)^n$ and the classifying $\infty$-category $\BN^r$, which will be used in the proof of Theorem \ref{thm: main infty replacement}. Here we will use $\varnothing$ to denote the initial object and $*$ to denote the final object in some $\infty$-category.

\blem\label{poset-colim-as-category}
    Let $P$ be a poset and $\cF$ be a family of its subsets, such that: 
    \benum
        \item\label{poset-colim-as-category-intersection-condition} $\cF$ is closed under finite intersection;
        \item\label{poset-colim-as-category-chain-condition} Every chain of finite length in $P$ is contained in some subset in $\cF$. 
    \eenu
    Then $P=\colim_{T\in\cF}T$ as $\infty$-categories, where we view posets as $\infty$-categories and view $\cF$ as a poset by subset inclusion. 
\elem

\bpf
    We prove the stronger statement that $P=\colim_{T\in\cF}T$ as simplicial animas. In other words, we check that for any map of posets $x\colon\Delta^n\to P$, the diagram
    $$\cF\ni T\mapsto\begin{cases}
        *,&x\text{ factors through }T;\\
        \varnothing,&\text{otherwise};
    \end{cases}$$
    has colimit $*$. As a diagram with values only $*$ and $\varnothing$, its colimit is the nerve of the part with value $*$, which is $\{T\in\cF\mid x(\Delta^n)\subseteq T\}$, as a subposet of $\cF$; now it is by (\ref{poset-colim-as-category-chain-condition}) nonempty and by (\ref{poset-colim-as-category-intersection-condition}) an inversely directed poset, so its nerve is contractible. 
\epf

\blem\label{directed-suspension-category}
    Let $C$ be an $\infty$-category, and let $C^\triangleleft$ and $C^\triangleright$ be the $\infty$-categories obtained from $C$ by adding an initial object $0$ and a final object $1$, respectively. In other words, $\Obj(C^\triangleleft)=\Obj(C)\sqcup\{0\}$, and the mapping animas are
    $$\Hom_{C^\triangleleft}(x,y)=\begin{cases}
        \Hom_C(x,y),&x,y\in C;\\
        *,&x=0;\\
        \varnothing,&\text{otherwise};
    \end{cases}$$
    similarly for $C^\triangleright$. Let $\Sigma C=C^\triangleleft\sqcup_CC^\triangleright$ be the pushout as $\infty$-categories. Then $\Obj(\Sigma C)=\Obj(C)\sqcup\{0,1\}$, and the mapping animas are
    $$\Hom_{\Sigma C}(x,y)=\begin{cases}
        \Hom_C(x,y),&x,y\in C;\\
        |C|,&x=0\text{ and }y=1;\\
        *,&x=0\text{ or }y=1\text{ but not both};\\
        \varnothing,&\text{otherwise};
    \end{cases}$$
    where $|C|$ means the nerve of the $\infty$-category $C$, i.e.\ the anima obtained from $C$ by inverting all its morphisms. 
\elem

\bpf
    This is a special case of \cite[Corollary 2.3]{ramzi-fun-with-pushouts}. 
\epf

\blem\label{pushout-diagram-product}
    Let $C$ be an $\infty$-category with pushouts and finite products, such that finite products distribute over pushouts. Let $T$ be a finite set, and for each $t\in T$ let
    $$\begin{tikzcd}
        x(0,t)\ar[r]\ar[d]&y(0,t)\ar[d]\\
        x(1,t)\ar[r]&y(1,t)
    \end{tikzcd}$$
    be a pushout square in $C$. Then the diagram $(\Delta^1)^T\times\Delta^1\to C$ given by
    $$(a,0)\mapsto\prod_{t\in T}x(a_t,t),\,(a,1)\mapsto\prod_{t\in T}y(a_t,t)$$
    is a colimit diagram. 
\elem

\bpf
    Note that if for each $t\in T$ we have a decomposition
    $$\begin{tikzcd}
        x(0,t)\ar[r]\ar[d]&z(0,t)\ar[r]\ar[d]&y(0,t)\ar[d]\\
        x(1,t)\ar[r]&z(1,t)\ar[r]&y(1,t)
    \end{tikzcd}$$
    where both small squares are pushouts, and we have proved the lemma for each of them in place of the large square, then we will have proved the lemma for the large square, because splicing two cubical colimit diagrams along a hyperface still gives a colimit diagram. Therefore, it suffices to treat the case where $x(0,t)=y(0,t)$ and $x(1,t)=y(1,t)$ for all but one $t\in T$, which is obvious because now the diagram is made up of the following faces which are pushouts by assumption,
    $$\begin{tikzcd}
        x(0,t)\times\prod_{u\ne t}x(a_u,u)\ar[r]\ar[d]&y(0,t)\times\prod_{u\ne t}x(a_u,u)\ar[d]\\
        x(1,t)\times\prod_{u\ne t}x(a_u,u)\ar[r]&y(1,t)\times\prod_{u\ne t}x(a_u,u)
    \end{tikzcd}$$
    where $a$ goes through $(\Delta^1)^{T\setminus t}$. 
\epf

\begin{example}\label{situation-cube}
    Let $S$ be a nonempty finite set and let $d=|S|$ be its order. Let $\cC_0$ be an $\infty$-category, $x=(x_0,x_1)\colon\partial\Delta^1\to\cC_0$ be two objects of it, and let $\cC=\cC_0\sqcup_{x,\partial\Delta^1}\Delta^1$ be the pushout in the $\infty$-category of $\infty$-categories, i.e.\ the $\infty$-category obtained from $\cC_0$ by adding a morphism from $x_0$ to $x_1$. We will consider the examples: 
    \benum
        \item\label{situation-example-cube} $\cC_0=\partial\Delta^1$ and $\cC=\Delta^1$;
        \item\label{situation-example-endomorphism} $\cC_0=*$ and $\cC=\BN$. 
    \eenu
\end{example}

\bprop\label{cube-boundary-category}
    Let $\partial(\Delta^1)^S$ be the boundary of the cube $(\Delta^1)^S$ viewed as an $\infty$-category. Rigorously speaking, 
    $$\partial(\Delta^1)^S=\colim_{T\subsetneq S}{(\Delta^1)^T\times(\partial\Delta^1)^{S\setminus T}}.$$
    Then the canonical maps $(\partial\Delta^1)^S\to\partial(\Delta^1)^S\to(\Delta^1)^S$ induce isomorphisms on animas of objects; namely, the anima of objects of $\partial(\Delta^1)^S$ is the set $\{0,1\}^S$. Abusing notation, let $0$ and $1$ denote the objects $0^S$ and $1^S$. Then for $x,y\in\partial(\Delta^1)^S$,
    $$\Hom_{\partial(\Delta^1)^S}(x,y)=\begin{cases}
        S^{d-2},&x=0\text{ and }y=1;\\
        \varnothing,&x\not\le y;\\
        *,&\text{otherwise}.
    \end{cases}$$
    Here $S^{d-2}$ means the $(d-2)$-dimensional sphere, and is more canonically the nerve of the $\infty$-category $(\Delta^1)^S\setminus\{0,1\}$. 
\eprop

\bpf
    Since $\partial\Delta^1=\{0,1\}$, we can break the colimit defining $\partial(\Delta^1)^S$ as
    $$\colim_{T\subsetneq S,\,i\in\{0,1\}^{S\setminus T}}(\Delta^1)^T,$$
    where $(T,i)\le(T',i')$ means $T\subseteq T'$ and $i\mapsto i'$ in the indexing poset, which we denote by $P$. Informally, $P$ is the inclusion poset of all faces of $(\Delta^1)^S$ of all dimensions, except the cube itself. Let $P_0$ and $P_1$ be subposets of $P$, consisting of faces away from $0$ and $1$, respectively. In other words,
    $$P_0=\{(T,i)\in P\mid i\ne0^{S\setminus T}\},\,P_1=\{(T,i)\in P\mid i\ne1^{S\setminus T}\}.$$
    Let $P_{01}=P_0\cap P_1$. Then Lemma \ref{poset-colim-as-category} obviously applies to the family $\{P_0,P_1,P_{01}\}$, giving $P=P_0\sqcup_{P_{01}}P_1$ as categories, and thus $\partial(\Delta^1)^S=\colim_{(T,i)\in P}(\Delta^1)^T$ is the pushout of $\colim_{(T,i)\in P_0}(\Delta^1)^T$ and $\colim_{(T,i)\in P_1}(\Delta^1)^T$ along $\colim_{(T,i)\in P_{01}}(\Delta^1)^T$. Now Lemma \ref{poset-colim-as-category} applies again, saying that the three colimits are the posets $(\Delta^1)^S\setminus0$, $(\Delta^1)^S\setminus1$, and $(\Delta^1)^S\setminus\{0,1\}$, respectively. This gives
    $$\partial(\Delta^1)^S=((\Delta^1)^S\setminus0)\sqcup_{(\Delta^1)^S\setminus\{0,1\}}((\Delta^1)^S\setminus1),$$
    and we can apply Lemma \ref{directed-suspension-category} to conclude. 
\epf

\bcor\label{category-cell}
    Let $E_S\subseteq\partial(\Delta^1)^S$ be the full sub-$\infty$-category containing only the two objects $0$ and $1$, so its mapping animas are $\Hom(0,0)=\Hom(1,1)=*$, $\Hom(0,1)=S^{d-2}$, and $\Hom(1,0)=\varnothing$. Then the following diagram is a pushout of $\infty$-categories,
    $$\begin{tikzcd}
        E_S\ar[r]\ar[d]&\partial(\Delta^1)^S\ar[d]\\
        \Delta^1\ar[r]&(\Delta^1)^S
    \end{tikzcd}$$
    where the top arrow is the inclusion functor, the right arrow is the canonical functor, and the other two arrows are the unique functors sending $0$ to $0$ and $1$ to $1$. 
\ecor

\bpf
    We prove the stronger statement that the above diagram is a pushout of simplicial animas. In other words, we check that the animas of functors
    $$\begin{tikzcd}
        \Fun(\Delta^n,E_S)\ar[r]\ar[d]&\Fun(\Delta^n,\partial(\Delta^1)^S)\ar[d]\\
        \Fun(\Delta^n,\Delta^1)\ar[r]&\Fun(\Delta^n,(\Delta^1)^S)
    \end{tikzcd}$$
    form pushout squares for all $n\in\N$. For this, we only need to check that for each $x\in\Fun(\Delta^n,(\Delta^1)^S)$ the following diagram is a pushout square of animas. 
    $$\begin{tikzcd}
        \Fun(\Delta^n,E_S)\times_{\Fun(\Delta^n,(\Delta^1)^S)}\{x\}\ar[r]\ar[d]&\Fun(\Delta^n,\partial(\Delta^1)^S)\times_{\Fun(\Delta^n,(\Delta^1)^S)}\{x\}\ar[d]\\
        \Fun(\Delta^n,\Delta^1)\times_{\Fun(\Delta^n,(\Delta^1)^S)}\{x\}\ar[r]&\{x\}
    \end{tikzcd}$$
    Now this is clear: if $x$ factors through the bottom arrow $\Delta^1\to(\Delta^1)^S$, then the horizontal arrows are isomorphisms, while if it does not, then the vertical arrows are isomorphisms. 
\epf

\begin{construction}
    For $T\subseteq S$, let $\cC^S_T=\cC^T\times\cC_0^{S\setminus T}$, so $\cC^S_\varnothing=\cC_0^S$ and $\cC^S_S=\cC^S$. The canonical map $\cC_0\to\cC$ makes $\cC^S_T$ functorial with respect to $T$, giving a functor from the poset of subsets of $S$ to the $\infty$-category of $\infty$-categories, i.e.\ a functor $(\Delta^1)^S\to\Cat$. For a downward-closed family $\cF$ of subsets of $S$, let 
    $$\cC^S_{\cF}=\colim_{T\in\cF}\cC^S_T,$$
    where $\cF$ is viewed as a full sub-$\infty$-category of $(\Delta^1)^S$ and the colimit is taken in $\Cat$. Then $\cC^S_{\{\varnothing\}}=\cC_0^S$ and $\cC^S_{\cP(S)}=\cC^S$, where $\cP(S)$ denotes the power set of $S$.  
\end{construction}

\bprop\label{cell-decomposition-power-category}
    Let $\cF_0$ be a downward-closed family of subsets of $S$, let $T$ be a nonempty subset of $S$ not in $\cF_0$ such that all its proper subsets are in $\cF_0$, and let $\cF=\cF_0\cup\{T\}$, so $\cF$ is also downward-closed. In this situation, we have the following pushout square of $\infty$-categories,
    $$\begin{tikzcd}
        E_T\times\cC_0^{S\setminus T}\ar[r]\ar[d]&\cC^S_{\cF_0}\ar[d]\\
        \Delta^1\times\cC_0^{S\setminus T}\ar[r]&\cC^S_{\cF}
    \end{tikzcd}$$
    where $E_T$ and the left arrow are as in Corollary \ref{category-cell}, the right arrow is the obvious map, and the horizontal arrows are given by the horizontal arrows in Corollary \ref{category-cell}, the functors $\partial\Delta^1\to\cC_0$ and $\Delta^1\to\cC$ in Example \ref{situation-cube}, and the definition of $\cC^S_{\cF}$. 
\eprop

\bpf
    Let $\cP=\cP(T)$ and $\cP_0=\cP(T)\setminus \{T\}$. Then $\cF=\cF_0\cup\cP$ and $\cF_0\cap\cP=\cP_0$, and Lemma \ref{poset-colim-as-category} applies, so $\cF=\cF_0\sqcup_{\cP_0}\cP$ as $\infty$-categories, and thus $\cC^S_{\cF}=\cC^S_{\cF_0}\sqcup_{\cC^S_{\cP_0}}\cC^S_{\cP}$. Therefore, to prove the proposition, it suffices to treat the case where $\cF=\cP$ and $\cF_0=\cP_0$. Canceling $\cC_0^{S\setminus T}$ and applying Corollary \ref{category-cell} with $T$ in place of $S$, we are reduced to showing that the following diagram is a pushout square. 
    $$\begin{tikzcd}
        \partial(\Delta^1)^T\ar[r]\ar[d]&\cC^T_{\cP_0}\ar[d]\\
        (\Delta^1)^T\ar[r]&\cC^T_{\cP}
    \end{tikzcd}$$
    Expanding the definitions of $\partial(\Delta^1)^T$ and $\cC^T_{\cP_0}$, we are reduced to showing that the diagram $\cP(T)\times\Delta^1\to\Cat$ given by
    $$(U,0)\mapsto(\Delta^1)^U\times(\partial\Delta^1)^{T\setminus U},\,(U,1)\mapsto\cC^U\times\cC_0^{T\setminus U}$$
    is a colimit diagram, which is Lemma \ref{pushout-diagram-product} with all the pushout squares as follows. 
    $$\begin{tikzcd}
        \partial\Delta^1\ar[r]\ar[d]&\cC_0\ar[d]\\
        \Delta^1\ar[r]&\cC
    \end{tikzcd}$$
\epf

\bcons\label{concrete-cell-unstable}
    We will only apply Proposition \ref{cell-decomposition-power-category} to the examples in Example \ref{situation-cube}, so we spell out here what Proposition \ref{cell-decomposition-power-category} implies in these cases: 
    \benum
        \item For any $\infty$-category $\cC$, giving a functor $X\colon(\Delta^1)^S\to\cC$ amounts to giving $2^{|S|}$ objects $X_\bullet\colon\{0,1\}^S\to\cC$, plus inductively on nonempty $T\subseteq S$, extending the maps $S^{|T|-2}\to\Hom(X_{0^T,i},X_{1^T,i})$ given by previous inductive steps along $S^{|T|-2}\to*$, where $i$ goes through elements of $\{0,1\}^{S\setminus T}$ so that $(0^T,i)$ and $(1^T,i)$ are elements of $\{0,1\}^S$. 
        \item For any $\infty$-category $\cC$, giving a functor $X\colon\BN^S\to\cC$ amounts to giving an object $X\in\cC$, plus inductively on nonempty $T\subseteq S$, extending the map $S^{|T|-2}\to\Hom(X,X)$ given by previous inductive steps along $S^{|T|-2}\to*$. 
    \eenu

    Write $d=|T|$. When $\cC$ is stable, the extension for $d=1$ amounts to choosing a morphism $X\to Y$. For $d\ge2$, a choice of basepoint gives a splitting $S^{d-2}\otimes X\simeq X\oplus\Sigma^{d-2}X$, so the extension amounts to choosing a nullhomotopy of the component $\Sigma^{d-2}X\to Y$.
\econs

\bpf[Proof of Theorem \ref{thm: main infty replacement}]
By Construction \ref{concrete-cell-unstable}(1) and \cite[Constructions 7.6 and 7.7]{Ye2026infty}, the extension data on the subcube with initial vertex $\ep$ and terminal vertex $\epp$ are represented by $D_{\ep,\epp}$, with boundary relation \eqref{eq: compatibility}. This gives the first assertion. For the second, apply Construction \ref{concrete-cell-unstable}(2) in the functor category of cubes. An action of $\BN^r$ on such a cube is precisely a choice of $r$ homotopy-coherently commutative operators on it.
\epf
% \bpf[Proof of Theorem \ref{thm: main infty replacement}]
%     The statement about the hypercube follows from Constructions \ref{cons: signs of cubes} and \ref{concrete-cell-unstable}, where the former provides a way to transfer a cube complex to a hypercube, and the latter provides the converse procedure. In particular, we use $D_{\ep,\epp}$ for $\ep={0^T,i},\epp={1^T,i}$ to extend the map $S^{|T|-2}\to\Hom(X_{0^T,i},X_{1^T,i})$ in the first case of Construction \ref{concrete-cell-unstable}.
    
%     The statement about the homotopy-coherently commutative operators follows from the same arguments, together with Construction \ref{cons: cube to BN}.
% \epf
\subsection{Periodic complexes}\label{subsec: Derived category of periodic chain complexes}

Periodic complexes are well-studied objects in homological algebra, but previous studies have focused only on their derived triangulated categories or dg-categories. In this subsection, we provide several equivalent definitions of the derived $\infty$-category of periodic complexes, beyond Definition \ref{defn: intro}. Some of these equivalences follow from \cite[\S 1.3]{HAlurie2017}, while others are $\infty$-categorical versions of results of Saito \cite{Saito2023cyclic}. The last description as a module category in Remark \ref{rem: periodic-complexes-as-dg-modules} relates the periodic and non-periodic settings and is useful for the construction of eigendecompositions. The advantage of the non-periodic derived $\infty$-category is that it has a natural t-structure.

We first recall the homotopy $\infty$-category and the derived $\infty$-category of non-periodic chain complexes. 

\bdefn[{\cite[Definition 1.3.2.1]{HAlurie2017}}]\label{def-derived}
    Let $\cA$ be an additive category. The \emph{homotopy $\infty$-category of $\cA$-complexes}, denoted $\K(\cA)$, is the dg-nerve $\Ndg_\bullet$ (\cite[Construction 1.3.1.6]{HAlurie2017}) of the dg-category $\Ch(\cA)$ of $\cA$-complexes, itself defined as follows: 
    \begin{itemize}
        \item An object $X\in\Ch(\cA)$ is of the form $(X_i,d_i)_{i\in\Z}$, where $X_i\in\cA$, $d_i\colon X_i\to X_{i-1}$, such that $d_{i-1}\circ d_i=0$ for all $i\in\Z$. 
        \item For $X,Y\in\Ch(\cA)$, the chain complex $\Hom(X,Y)_\bullet$ is defined as
        $$\Hom(X,Y)_p=\prod_{i\in\Z}\Hom_\cA(X_i,Y_{i+p});\;(df)(x)=d(f(x))-(-1)^pf(dx).$$
    \end{itemize}
    We know that $\K(\cA)$ is a stable $\infty$-category by \cite[Proposition 1.3.2.10]{HAlurie2017}. 

    If $\cA$ is an abelian category, then one can consider the full sub-$\infty$-category $\Ac(\cA)\subset\K(\cA)$ of acyclic complexes, i.e.\ those with trivial homology. The \emph{derived $\infty$-category} $\De(\cA)$ of $\cA$ is then defined as the quotient $\K(\cA)/\Ac(\cA)$, i.e.\ the initial stable $\infty$-category under $\K(\cA)$ such that $\Ac(\cA)$ is mapped to zero. 
\edefn

Now we adapt the above definition to periodic complexes. 

\bdefn\label{def-periodic}
    Let $m$ be a positive integer and $\cA$ be an additive category. The \emph{homotopy $\infty$-category of $m$-periodic $\cA$-complexes}, denoted $\K_m(\cA)$, is the dg-nerve of the dg-category $\Ch_m(\cA)$ of $m$-periodic $\cA$-complexes, itself defined as follows: 
    \begin{itemize}
        \item An object $X\in\Ch_m(\cA)$ is of the form $(X_i,d_i)_{i\in\Z/m}$, where $X_i\in\cA$, $d_i\colon X_i\to X_{i-1}$, such that $d_{i-1}\circ d_i=0$ for all $i\in\Z/m$. 
        \item For $X,Y\in\Ch_m(\cA)$, the chain complex $\Hom(X,Y)_\bullet$ is defined as
        $$\Hom(X,Y)_p=\bigoplus_{i\in\Z/m}\Hom_{\cA}(X_i,Y_{i+p});\;(df)(x)=d(f(x))-(-1)^pf(dx).$$
        Note that $p$ is not an element of $\Z/m$ but an integer, so $(-1)^p$ makes sense. 
    \end{itemize}
    It is easy to see that the proof of \cite[Proposition 1.3.2.10]{HAlurie2017} carries over and shows that $\K_m(\cA)$ is a stable $\infty$-category. 

    If $\cA$ is an abelian category, one can likewise consider the full sub-$\infty$-category $\Ac_m(\cA)\subset\K_m(\cA)$ of acyclic complexes, i.e.\ those with trivial homology. The \emph{$m$-periodic derived $\infty$-category} $\De_m(\cA)$ is then defined as the quotient $\K_m(\cA)/\Ac_m(\cA)$, i.e.\ the initial stable $\infty$-category under $\K_m(\cA)$ such that $\Ac_m(\cA)$ is mapped to zero. 
\edefn

\brem
    In Definitions \ref{def-derived} and \ref{def-periodic}, taking quotients by acyclic complexes is the same as inverting quasi-isomorphisms, because quasi-isomorphisms are exactly those morphisms with acyclic cofibers, and in a stable $\infty$-category, a morphism is an isomorphism if and only if its cofiber is zero. Hence the definitions are equivalent to the localizations $\K(\cA)[W^{-1}]$ and $\K_m(\cA)[W_m^{-1}]$ from \cite[Definition 1.3.4.1]{HAlurie2017}, respectively, where $W$ and $W_m$ are the collections of quasi-isomorphisms.
\erem

\brem
    Keep the notation of Definitions \ref{def-derived} and \ref{def-periodic}. By definition of dg-nerves, a map $f\colon X\to Y$ in $\K(\cA)$ (resp.\ $\K_m(\cA)$) is represented by a $0$-cycle of the chain complex $\Hom(X,Y)_\bullet$, i.e.\ $f\in\Hom(X,Y)_0$ with $df=0$. It is easy to verify that this corresponds to a chain map $f\colon X\to Y$ in $\Ch(\cA)$ (resp.\ $\Ch_m(\cA)$). Also easy is that an equality between two maps in $\K(\cA)$ represented by chain maps $f,g\colon X\to Y$ is represented by a chain homotopy between $f$ and $g$, an equality between equalities is represented by homotopy between homotopies, and so on (cf.\ \cite[Example 1.3.1.8, Remark 1.3.2.2]{HAlurie2017}). 
\erem

The advantage of Definition \ref{def-periodic} over Definition \ref{defn: intro} is that the former fits better into the framework of higher algebra:

\begin{notation}\label{not: chain complex}
    Given a commutative ring $\bK$ with identity, we write $\Mod(\bK)$ for the abelian category of (possibly infinitely presented) $\bK$-modules and write $\De(\bK)=\De(\Mod(\bK))$. The same abbreviation applies to $\De_m(\bK)$ and a \emph{$\bK$-complex} means a \emph{$\Mod(\bK)$}-complex.
    
    % As we are particularly interested in $\bK=\C$ and $m=2$, we will write $\De_2=\De_2(\C)$ in subsequent sections. A $\C$-complex is sometimes simply called a (chain) complex.
\end{notation}

%If you think that we need to elaborate on the symmetric monoidal structure of the dg-category Ch(K) and of the localization functor Ch(K)→D(K), I can add a remark here. 

\brem[periodic complexes as dg-modules]\label{rem: periodic-complexes-as-dg-modules}
    Let $\bK$ be a commutative ring with identity and $m$ be a positive integer. Consider the graded Laurent algebra $\bK[\beta^\pm]$ with $\beta$ in degree $m$, viewed as a differential graded algebra with zero differentials. Then it is obvious that $m$-periodic complexes are the same as right dg-modules of the dg-algebra $\bK[\beta^\pm]$, i.e.\ 
    $$\Ch_m(\bK)=\RMod_{\bK[\beta^\pm]}(\Ch(\bK)).$$
    Applying \cite[Theorem 4.3.3.17]{HAlurie2017} with $\mathbf{A}=\Ch(\bK)$, $A=\bK$ and $B=\bK[\beta^\pm]$, we see that
    $$\De_m(\bK)=\RMod_{\bK[\beta^\pm]}(\De(\bK)),$$
    i.e.\ $\De_m(\bK)$ is equivalent to the $\infty$-category of right $\bK[\beta^\pm]$-modules in the $\infty$-category $\De(\bK)$. 

    When $m$ is even or $2=0$ in $\bK$, $\bK[\beta^\pm]$ is a commutative dg-algebra. In this case, there is no difference between left and right dg-modules, and moreover $\De_m(\bK)=\Mod_{\bK[\beta^\pm]}(\De(\bK))$ is a symmetric monoidal $\infty$-category by \cite[Theorem 4.5.2.1]{HAlurie2017}, using the derived tensor product on $\De(\bK)$ from \cite[Remark 7.1.2.12]{HAlurie2017}.
\erem

In the rest of this subsection, we show that the derived $\infty$-categories in Definitions \ref{defn: intro} and \ref{def-periodic} are equivalent (Theorem \ref{projective-resolution-equivalence}) under the assumption in Notation \ref{not: assumption}. Moreover, this result also implies the symmetric monoidal structure defined in Remark \ref{rem: periodic-complexes-as-dg-modules} is equivalent to that in \cite{Ye2026infty}.

\bdefn[{cf.\ the dual of \cite[\texttt{070H}]{stacks}}]\label{k-projective}
    Let $m$ be a positive integer and $\cA$ an abelian category. An object $P$ in $\K(\cA)$ (resp.\ $\K_m(\cA)$) is called \emph{projective}, if for all $X$ in $\Ac(\cA)$ (resp.\ $\Ac_m(\cA)$), we have $\Hom_{\K(\cA)}(P,X)=0$ (resp.\ $\Hom_{\K_m(\cA)}(P,X)=0$). Note that projective objects form a stable sub-$\infty$-category of $\K(\cA)$ (resp.\ $\K_m(\cA)$) that is closed under colimits.
\edefn

\bprop\label{hom-from-projective-in-derived}
    Keep the notation of Definition \ref{k-projective}, so $P$ is a projective ($m$-periodic) complex. Then for all $Y$ in $\K(\cA)$ (resp.\ $\K_m(\cA)$), we have \[\Hom_{\De(\cA)}(P,Y)=\Hom_{\K(\cA)}(P,Y)~(\text{resp.}\ \Hom_{\De_m(\cA)}(P,Y)=\Hom_{\K_m(\cA)}(P,Y)).\]
\eprop

\bpf
    This is a general fact: for a stable $\infty$-category $\cC$ and its stable full sub-$\infty$-category $\cD$, if $P\in\cC$ satisfies $\Hom(P,X)=0$ for all $X\in\cD$, then we have $\Hom_{\cC}(P,Y)=\Hom_{\cC/\cD}(P,Y)$ for all $Y\in\cC$. 
\epf

\bcor
    Keep the notation of Definition \ref{k-projective}, and let $f\colon P\to Q$ be a map in $\K(\cA)$ (resp.\ $\K_m(\cA)$) where both $P$ and $Q$ are projective. If $f$ is an isomorphism in $\De(\cA)$ (resp.\ $\De_m(\cA)$), then it is an isomorphism in $\K(\cA)$ (resp.\ $\K_m(\cA)$), i.e.\ a chain homotopy equivalence. 
\ecor

Recall that an object $P$ of an abelian category $\cA$ is said to be \emph{projective} if $\Hom_{\cA}(P,-)$ is an exact functor. 

\bprop\label{projectivity-bounded-right}
    Let $\cA$ be an abelian category. If an $\cA$-complex $P$ has all its terms projective and is bounded on the right (i.e.\ $P_i=0$ for $i\ll0$), then $P$ is projective in $\K(\cA)$. 
\eprop

\bpf
    Let $X$ be an acyclic $\cA$-complex. We need to show that $\Hom_{\K(\cA)}(P,X)=0$, i.e.\ every chain map $P\to X$ is nullhomotopic. Let $f=(f_i\colon P_i\to X_i)_{i\in\Z}$ be a chain map. We inductively construct the nullhomotopy $h$ of $f$ as follows: 
    \begin{itemize}
        \item For $i\ll0$, since $P_i=0$, we have no choice but to set $h_i=0$. 
        \item Suppose that we have constructed $h_j$ that satisfy $d_{j+1}h_j+h_{j-1}d_j=f_j$ for $j<i$. We want to do the same thing for $i$. Namely, we want to find $h_i\colon P_i\to X_{i+1}$ such that $d_{i+1}h_i=f_i-h_{i-1}d_i$. Since $P_i$ is projective, it suffices to see that the image of $f_i-h_{i-1}d_i$ is contained in that of $d_{i+1}$; since $X$ is acyclic, this is equivalent to saying that $d_i(f_i-h_{i-1}d_i)=0$, which can be seen quite easily by the inductive hypothesis. \qedhere
    \end{itemize}
\epf

\bprop\label{homotopy-complex-colimit-termwise}
    Let $\cA$ be an abelian category. Let
    $$P(0)\to P(1)\to\cdots$$
    be an $\N$-diagram of $\cA$-complexes, i.e.\  a functor from the poset $(\N,\le)$ to the category $\Ch(\cA)$. Assume that for every $i\in\Z$, the diagram on the $i$-th term
    $$P(0)_i\to P(1)_i\to\cdots$$
    is eventually constant, i.e.\ there exists an $N$ (depending on $i$) such that
    $$P(N)_i\xrightarrow{\sim}P(N+1)_i\xrightarrow{\sim}\cdots.$$
    Denote this $P(N)_i$ by $P_i$, and denote the obvious chain complex $(P_i)_{i\in\Z}$ by $P$. Then
    $P=\colim_{n\to\infty}P(n)$ in $\K(\cA)$, i.e.\ the colimit of $P(n)$ in $\K(\cA)$ can be computed termwise. 
\eprop

\bpf
    By definition of colimits, we need to show that for all $Y\in\K(\cA)$, 
    $$\Hom_{\K(\cA)}(P,Y)=\lim_{n\to\infty}\Hom_{\K(\cA)}(P(n),Y),$$
    as limits of animas. By definition of dg-nerves, $\Hom_{\K(\cA)}(X,Y)$ can be computed as the $\Omega^\infty$ of the complex $\Hom(X,Y)_\bullet\in\De(\Z)$, and since $\Omega^\infty$ preserves limits, it suffices to show the above equality in $\De(\Z)$. Now the proposition follows from Lemma \ref{derived-complex-limit-termwise}, since $\limone$ of eventually constant diagrams are zero. 
\epf

\blem\label{derived-complex-limit-termwise}
    Let
    $$A(0)\gets A(1)\gets\cdots$$
    be an $\N^\opp$-diagram of chain complexes of abelian groups. Assume that for all $i\in\Z$, 
    $$\limone_{n\to\infty}A(n)_i=0.$$
    Then the limit $\lim_{n\to\infty}A(n)$ in $\De(\Z)$ can be computed as the termwise limit. 
\elem

\bpf
    Here $\limone$ means the $R^1\lim$ in \cite[\texttt{091D}]{stacks}. The lemma follows immediately from \cite[\texttt{07KC}, \texttt{091D}, \texttt{0CQD}]{stacks}. 
\epf

\bdefn\label{defn: global dimension}
An abelian category $\cA$ has \emph{enough projectives} if every object of $\cA$ admits an epimorphism from a projective object. Its \emph{global dimension} is
\[
\gldim\cA=\sup\{i\in\N\mid \Ext^i(X,Y)\ne0\text{ for some }X,Y\in\cA\}.
\]
If $\cA$ has enough projectives, then $\gldim\cA\le d$ if and only if every object of $\cA$ admits a projective resolution of length at most $d$. For a commutative ring $\bK$, we write $\gldim\bK=\gldim\Mod(\bK)$.
\edefn
\begin{notation}\label{not: assumption}
    Until the end of this subsection, let $\cA$ be an abelian category with enough projectives and of finite global dimension, say $\gldim\cA=d\in\N$. Let $\Pj$ denote the full subcategory of its projective objects.
\end{notation} 

\blem\label{nonperiodic-projective-resolution}
    For every $\cA$-complex $X$, there exists an $\cA$-complex $P$ and a quasi-isomorphism $P\to X$, such that $P$ is projective in $\K(\cA)$, and all the terms of $P$ are projective. In particular, the natural functor $\K(\Pj)\to\De(\cA)$ is essentially surjective. 
\elem

\bpf
    Consider the sequence of natural truncations of $X$
    $$\cdots\to\tau_{\ge1}X\to\tau_{\ge0}X\to\tau_{\ge-1}X\to\cdots$$
    where $\tau_{\ge n}X$ is defined as the complex 
    $$\cdots\to X_{n+2}\to X_{n+1}\to Z_n(X)\to0\to0\to\cdots$$
    which has homologies of degree $\ge n$ the same as those of $X$ and other homologies zero. Now it suffices to construct a sequence $(P(n))_{n\in\N}$ of $\cA$-complexes that are projective in $\K(\cA)$, have projective terms, and fit into a diagram of complexes
    $$\begin{tikzcd}
        P(0)\ar[r]\ar[d]&P(1)\ar[r]\ar[d]&\cdots\\
        \tau_{\ge0}X\ar[r]&\tau_{\ge-1}X\ar[r]&\cdots
    \end{tikzcd}$$
    commuting in $\K(\cA)$, whose upper row satisfies the assumption of Proposition \ref{homotopy-complex-colimit-termwise} and whose vertical maps are quasi-isomorphisms, because then the termwise colimit of complexes $P=\colim_{n\to\infty}P(n)$ computes the colimit in $\K(\cA)$, giving a map $P\to X$ that meets all our requirements. 

    We do the construction inductively as follows:
    \begin{description}
        \item[Base Case] This is the well-known projective resolution of right bounded complexes. See for example \cite[\texttt{0646}]{stacks}. 
        \item[Induction] Assume that we have constructed the desired $P(n)\to\tau_{\ge-n}X$. Note that in $\De(\cA)$ we have $\fib(\tau_{\ge-n}X\to\tau_{\ge-n-1}X)=H_{-n-1}(X)[-n-2]$, by the homology long exact sequence. By our assumption that $\cA$ has global dimension $d$, we can find a quasi-isomorphism $Q\to H_{-n-1}(X)[-n-2]$ such that terms of $Q$ are projective and are nonzero only in the range $[-n-2,d-n-2]$. Since $Q\cong H_{-n-1}(X)[-n-2]$ and $P(n)\cong\tau_{\ge-n}X$ in $\De(\cA)$, the natural map $H_{-n-1}(X)[-n-2]\to\tau_{\ge-n}X$ can be viewed as a map $Q\to P(n)$ in $\De(\cA)$, which by Propositions \ref{hom-from-projective-in-derived} and \ref{projectivity-bounded-right} is also a map $Q\to P(n)$ in $\K(\cA)$, and thus can be represented by a chain map of complexes. Now let $P(n+1)$ be the mapping cone of it. Then by construction, the chain map $P(n)\to P(n+1)$ is an isomorphism on degrees $>d-n-1$, and we have a diagram
        $$\begin{tikzcd}
            P(n)\ar[r]\ar[d]&P(n+1)\ar[d]\\
            \tau_{\ge-n}X\ar[r]&\tau_{\ge-n-1}X
        \end{tikzcd}$$
        commuting in $\De(\cA)$. By Propositions  \ref{hom-from-projective-in-derived} and \ref{projectivity-bounded-right} again, $P(n+1)$ is projective in $\K(\cA)$, the vertical maps above are also maps in $\K(\cA)$, and the diagram commutes in $\K(\cA)$, so our inductive construction is done. \qedhere
    \end{description}
\epf

\blem\label{periodic-projective-resolution}
    For every $m$-periodic $\cA$-complex $X$, there exists an $m$-periodic $\cA$-complex $P$ and a quasi-isomorphism $P\to X$, such that $P$ is projective in $\K_m(\cA)$, and all the terms of $P$ are projective. In particular, the natural functor $\K_m(\Pj)\to\De_m(\cA)$ is essentially surjective. 
\elem

\bpf
    For $Q\in\cA$, let $Q[0]$ denote the complex with $Q$ in degree $0$, and with all differentials zero. We first show that if $Q\in\cA$ is projective, then $Q[0]\in\K_m(\cA)$ is projective. Take an acyclic $X\in\K_m(\cA)$ and a chain map $f\colon Q[0]\to X$. We want to show that $f$ is nullhomotopic. Note that $f$ is essentially a map $f_0\colon Q\to X_0$ with $d_0f_0=0$. By acyclicity of $X$ and projectivity of $Q$, this implies that $f_0$ factors through $d_1\colon X_1\to X_0$, giving the desired nullhomotopy. 

    The acyclic case is immediate. We do induction on the length of an interval containing the nonzero homology degrees of $X$.
    % We do induction on nonzero homologies of $X$ as follows:
    \begin{description}
        \item[Base Case] Assume that $X$ has only one nonzero homology. In this case we do induction on the projective dimension of that homology. Without loss of generality, say $A=H_0(X)\ne0$. If $A$ has projective dimension $0$, i.e.\ $A$ is projective, then there will be a section $A\to Z_0(X)$ to the quotient map $Z_0(X)\to Z_0(X)/B_0(X)=H_0(X)$, giving rise to a quasi-isomorphism $A[0]\to X$. If $A$ has positive projective dimension, then take a projective $Q$ and a surjection $Q\twoheadrightarrow A$. Similarly, there is a map $Q\to Z_0(X)$ lifting $Q\twoheadrightarrow A=Z_0(X)/B_0(X)$, giving rise to a chain map $Q[0]\to X$. It is easy to see that $\fib(Q[0]\to X)$ still has only one nonzero homology, whose projective dimension is less than that of $A$, so by the induction hypothesis there is a projective complex $P$ with all its terms projective, and a quasi-isomorphism $P\to\fib(Q[0]\to X)$. This gives a quasi-isomorphism $\cofib(P\to Q[0])\to X$. Since $\cofib(P\to Q[0])$ can be represented by the mapping cone which has all the terms projective, and it is projective because projective objects form a stable sub-$\infty$-category, we are done. 
        \item[Induction] Assume that $X$ has nonzero homology in more than one degree. Without loss of generality, say $A=H_0(X)\ne0$. Take a positive integer $n<m$ such that nonzero homologies of $X$ are concentrated on $0,1,\ldots,n$ and do induction on $n$. Take a projective $Q$ and a surjection $Q\twoheadrightarrow A$. By projectivity, there is a map $Q\to Z_0(X)$ lifting $Q\twoheadrightarrow A=Z_0(X)/B_0(X)$, giving rise to a chain map $Q[0]\to X$. By the long exact sequence of homologies, it is easy to see that nonzero homologies of $\fib(Q[0]\to X)$ are concentrated on $0,1,\ldots,n-1$. Now using the induction hypothesis, the same reasoning as in the Base Case works, and we are done. \qedhere
    \end{description}
\epf

\begin{theorem}\label{projective-resolution-equivalence}
    Let $m$ be a positive integer. Let $P$ be an ($m$-periodic) $\cA$-complex with all its terms projective. Then $P$ is projective in $\K(\cA)$ (resp.\ $\K_m(\cA)$). The natural functor $\K(\Pj)\to\De(\cA)$ (resp.\ $\K_m(\Pj)\to\De_m(\cA)$) is an equivalence. 
\end{theorem}

\bpf
    We first deal with the non-periodic case. By Lemma \ref{nonperiodic-projective-resolution}, there exists an $\cA$-complex $P'$ which is projective in $\K(\cA)$ and has all its terms projective, and a quasi-isomorphism $P'\to P$. So projectivity of $P$ is equivalent to that of the mapping cone of $P'\to P$, which is an acyclic complex with all the terms projective. Therefore, we can assume that $P$ is acyclic, without loss of generality. In this case we claim that $P=0$ in $\K(\cA)$, i.e.\ $P$ is homotopic to the zero complex. For this it suffices to show that the short exact sequence
    $$0\to Z_i(P)\to P_i\to B_{i-1}(P)=Z_{i-1}(P)\to0$$
    splits for every $i$, because then $P$ will be the direct sum of shifts of two-term complexes $0\to Z_i(P)\xrightarrow{\id}Z_i(P)\to0$, which themselves are clearly nullhomotopic. To that end, consider the functor $\Ext^n(Z_i(P),-)$ for $n>0$. By the long exact sequence of $\Ext$ associated to the above short exact sequence, we see that
    $$\Ext^n(Z_i(P),-)=\Ext^{n+1}(Z_{i-1}(P),-)=\Ext^{n+2}(Z_{i-2}(P),-)=\cdots$$
    Since $\cA$ is of finite global dimension and the equalities continue infinitely, all the $\Ext$ above have to be zero. This implies that $Z_i(P)$ is projective for every $i$ and hence the above short exact sequence splits. Now $\K(\Pj)\to\De(\cA)$ is fully faithful due to Proposition \ref{hom-from-projective-in-derived}, and essentially surjective due to Lemma \ref{nonperiodic-projective-resolution}, so it is an equivalence. 

    The proof for the periodic case is the same, with Lemma \ref{periodic-projective-resolution} in place of Lemma \ref{nonperiodic-projective-resolution}. 
\epf

\bcor\label{cor: enough proj}
    Let $m$ be a positive integer. Let $\cA_0\subseteq\cA$ be a full abelian subcategory where $\Pj_0=\Pj\cap\cA_0$ is a family of enough projectives in $\cA_0$, i.e.\ for all $A\in\cA_0$ there exists $Q\in\Pj_0$ with a surjection $Q\twoheadrightarrow A$. Then the natural functor $\De(\cA_0)\to\De(\cA)$ (resp.\ $\De_m(\cA_0)\to\De_m(\cA)$) is fully faithful. 
\ecor

\bpf
    By Theorem \ref{projective-resolution-equivalence} this reduces to full faithfulness of $\K(\Pj_0)\to\K(\Pj)$ (resp.\ $\K_m(\Pj_0)\to\K_m(\Pj)$), which is obvious by definition of dg-nerves. 
\epf

\bexmp\label{exmp: periodic-derived-compact}
    In Corollary \ref{cor: enough proj}, $\cA$ can be the ordinary module category $\Mod(\bK)$ of a commutative Noetherian ring $\bK$ of finite global dimension. The full subcategory $\cA_0=\Mod(\bK)^c$ of compact objects is an abelian subcategory. Then $\De_m(\cA)$ is compactly generated, and the natural functor $\De_m(\cA_0)\to\De_m(\cA)$ is fully faithful with essential image exactly the compact objects of $\De_m(\cA)$. Note that we cannot replace the periodic assumption because the compact objects in $\De(\cA)$ should be bounded, while the essential image of $\De(\cA_0)$ contains unbounded complexes.
\eexmp

\brem\label{rem: same sym mono}
    Suppose $\bK$ has finite global dimension and that $m$ is even or $2=0$ in $\bK$. For $\cA=\Mod(\bK)$, the equivalence in Theorem \ref{projective-resolution-equivalence} identifies the symmetric monoidal structure of \cite[Construction 3.54]{Ye2026infty} with that of Remark \ref{rem: periodic-complexes-as-dg-modules}. Indeed, both are symmetric monoidal localizations at quasi-isomorphisms of the ordinary category of degreewise projective $m$-periodic complexes with its Koszul tensor product, which computes the derived tensor product over $\bK[\beta^\pm]$ \cite[Proposition 4.1.7.4 and Remark 4.1.7.5]{HAlurie2017}.
\erem

\subsection{Eigendecomposition}\label{subsec: Eigenspace decomposition}

Let $\bK$ be an algebraically closed field. In practice, our periodic complexes will be acted on by $r$ homotopy-coherently commutative operators, i.e.\ they will be linear over $A_r=\bK[T_1,\ldots,T_r]$ as in Definition \ref{defn: intro}, and we would like to take the eigendecomposition with respect to these operators. In other words, viewing geometrically rather than algebraically, our complexes are quasicoherent complexes on $\Ga^r=\Spec(A_r)$, the $r$-fold product of the additive group, with the periodicity action of Remark \ref{rem: periodic-complexes-as-dg-modules}. Recall by \cite[Theorem 7.1.2.13]{HAlurie2017} and \cite[Proposition 2.2.3.3]{SAGlurie2018} the equivalence $\De(A_r)\simeq\QCoh(\Ga^r)$ of $\infty$-categories. In this subsection we give decompositions of certain quasicoherent complexes according to their supports, symmetric monoidally with respect to the convolution product, which should be well-known among experts but we cannot find a reference. 

% The existence of the decomposition in Theorem \ref{thm: main eigenspace decomposition} follows from finite support and \cite[\texttt{0AEF}]{stacks}, as explained below. Our main task is to establish its compatibility with the symmetric monoidal structures, including in the periodic case.
The existence of the eigendecomposition in Theorem \ref{thm: main eigenspace decomposition} follows from finite support and \cite[\texttt{0AEF}]{stacks}, as explained below. Our main task is to establish its compatibility with the symmetric monoidal structures, including in the periodic case.

In this subsection, we develop the theory of such decompositions in a more general setting, namely with $\Ga^r$ replaced by a general commutative algebraic group $G$ over $\bK$, which we fix from now throughout the subsection.

\brem
    By \cite[Proposition 9.6.1.1]{SAGlurie2018}, for a general qcqs (i.e.\ quasicompact and quasiseparated) scheme $X$, the $\infty$-category $\QCoh(X)$ of quasicoherent $\cO_X$-complexes is compactly generated. For a morphism $f\colon X\to Y$ of such schemes, there is an adjunction
    $$f^*\colon\QCoh(Y)\rightleftarrows\QCoh(X):\!f_*$$
    so $f^*$ preserves colimits and $f_*$ preserves limits. When $f$ is a qcqs morphism, $f_*$ also preserves colimits by \cite[Corollary 3.4.2.2]{SAGlurie2018}. Beware that our $(f^*,f_*)$ means the derived pullback and pushforward; the underived functors do not in general preserve quasi-isomorphisms. 
\erem

In this subsection, we will mainly use the pushforward functoriality of $\QCoh$, i.e.\ view $\QCoh$ as a covariant functor from qcqs schemes to compactly generated $\infty$-categories.

% In this subsection, we will mainly use the pushforward functoriality of $\QCoh$, i.e.\ view $\QCoh$ as a covariant functor from qcqs schemes to compactly generated presentable $\infty$-categories.
\bcons[convolution on $\QCoh(G)$]
%     By \cite[Corollary 9.4.2.3]{SAGlurie2018}, the functor $\QCoh$ from finite type $\bK$-schemes to $\Mod(\bK)$-linear presentable $\infty$-categories is symmetric monoidal. Therefore, it sends our commutative algebraic group $G$ over $\bK$ to a commutative group object in $\Mod(\bK)$-linear presentable $\infty$-categories. This gives $\QCoh(G)$ a symmetric monoidal structure with antipode, which we call the \emph{convolution product}.
    Let $a\colon G\times G\to G$ be addition and $0\colon\Spec(\bK)\to G$ the identity. The \emph{convolution product} is $A\otimes B=a_*(A\boxtimes B)$, where $\boxtimes$ is the derived external tensor product. By the group law and \cite[Corollary 9.4.2.3]{SAGlurie2018}, this gives $\QCoh(G)$ a symmetric monoidal structure with unit $0_*\bK$. 
\econs

In this generality, ``complexes that decompose into their generalized eigenspaces'' are quasicoherent complexes with $0$-dimensional support. Before studying them we first recall some facts on coidempotents, appearing for example in \cite[\S 2]{aoki-2023-sheaves-spectrum-adjunction}: 

\brem\label{coidem=coloc}
    Let $(\cC,\otimes,\mathbf{1})$ be a symmetric monoidal $\infty$-category. The following notions are equivalent: 
    \begin{description}
        \item[Coidempotent] An object $c\colon C\to\mathbf{1}$ in the overcategory $\cC_{/\mathbf{1}}$ such that $C\otimes c\colon C\otimes C\to C\otimes\mathbf{1}=C$ is an isomorphism. 
        \item[Coidempotent coalgebra] A cocommutative coalgebra $C$ in $\cC$ whose comultiplication $C\to C\otimes C$ is an isomorphism.
%         \item[Smashing colocalization] A symmetric monoidal full sub-$\infty$-category $\cC'\subseteq\cC$ whose inclusion functor has a right adjoint $R\colon\cC\to\cC'$, for which the natural transformation $R(\mathbf{1})\otimes-\to R(-)$ is an isomorphism.
        \item[Smashing colocalization] A full sub-$\infty$-category $\cC'\subseteq\cC$ closed under tensoring with objects of $\cC$, whose inclusion functor has a right adjoint $R\colon\cC\to\cC'$, for which the natural transformation $R(\mathbf{1})\otimes-\to R(-)$ is an isomorphism. 
    \end{description}
%     The equivalence goes as follows: the coidempotent $c\colon C\to\mathbf{1}$ corresponds to the full sub-$\infty$-category of objects $X\in\cC$ such that $X\otimes c\colon X\otimes C\to X$ is an isomorphism, whose inclusion functor has the right adjoint $-\otimes C$; conversely, a smashing localization $\cC'\subseteq\cC$ with right adjoint $R$ corresponds to the coidempotent $R(\mathbf{1})$.
    The equivalence goes as follows: the coidempotent $c\colon C\to\mathbf{1}$ corresponds to the full sub-$\infty$-category of objects $X\in\cC$ such that $X\otimes c\colon X\otimes C\to X$ is an isomorphism, whose inclusion functor has the right adjoint $-\otimes C$; conversely, a smashing colocalization $\cC'\subseteq\cC$ with right adjoint $R$ corresponds to the coidempotent $R(\mathbf{1})$. 
\erem

\bexmp\label{local-cohomology-multiple}
%     Recall from \cite[Definition 7.1.5.1]{SAGlurie2018} the full subcategory $\QCoh_Z(G)\subseteq\QCoh(G)$ for a closed subset $Z\subseteq G$, and from \cite[Proposition 7.2.3.1]{SAGlurie2018} that the inclusion has a right adjoint $\RG_Z$ called \emph{local cohomology} (cf.\ \cite[\texttt{0952}]{stacks}), which is a smashing localization corresponding to the coidempotent $\RG_Z(\cO_G)\in\QCoh_Z(G)\subseteq\QCoh(G)$.
    Recall from \cite[Definition 7.1.5.1]{SAGlurie2018} the full subcategory $\QCoh_Z(G)\subseteq\QCoh(G)$ for a closed subset $Z\subseteq G$, and from \cite[Proposition 7.2.3.1]{SAGlurie2018} that the inclusion has a right adjoint $\RG_Z$ called \emph{local cohomology} (cf.\ \cite[\texttt{0952}]{stacks}), which is a smashing colocalization for the usual derived tensor product over $\cO_G$, corresponding to the coidempotent $\RG_Z(\cO_G)\in\QCoh_Z(G)\subseteq\QCoh(G)$. 

%     For disjoint closed $Z,Z'\subseteq G$, $\RG_Z(\cO_G)\otimes\RG_{Z'}(\cO_G)=0$. Therefore, for a family $(Z_\lambda)_{\lambda\in\Lambda}$ of disjoint closed subsets of $G$, the direct sum $\bigoplus_{\lambda\in\Lambda}\RG_{Z_\lambda}(\cO_G)\in\QCoh(G)$ is also a coidempotent, and corresponds to a smashing localization, which we denote by $\RG_{(Z_\lambda)_{\lambda\in\Lambda}}$ and whose image we denote by $\QCoh_{(Z_\lambda)_{\lambda\in\Lambda}}(G)$.
    For disjoint closed $Z,Z'\subseteq G$, $\RG_Z(\cO_G)\otimes\RG_{Z'}(\cO_G)=0$. Therefore, for a family $(Z_\lambda)_{\lambda\in\Lambda}$ of disjoint closed subsets of $G$, the direct sum $\bigoplus_{\lambda\in\Lambda}\RG_{Z_\lambda}(\cO_G)\in\QCoh(G)$ is also a coidempotent, and corresponds to a smashing colocalization, which we denote by $\RG_{(Z_\lambda)_{\lambda\in\Lambda}}$ and whose image we denote by $\QCoh_{(Z_\lambda)_{\lambda\in\Lambda}}(G)$. 
\eexmp

\bdefn[$0$-dimensionally supported complexes]
%     In the above Example \ref{local-cohomology-multiple}, take $(Z_\lambda)_{\lambda\in\Lambda}$ to be all the closed points of the algebraic group $G$, which are also all the $\bK$-points of $G$ as $\bK$ is algebraically closed. We denote by $\RGzd$ the resulting smashing localization and $\QCohzd(G)$ its image. By definition, $A=\bigoplus_{\lambda\in G(\bK)}\RG_{\lambda}(A)$ for all $A\in\QCohzd(G)$, where $\RG_{\lambda}(A)$ is supported on the closed point $\lambda$. Therefore, $\QCohzd(G)=\prod_{\lambda\in G(\bK)}\QCoh_{\lambda}(G)$.
    In the above Example \ref{local-cohomology-multiple}, take $(Z_\lambda)_{\lambda\in\Lambda}$ to be all the closed points of the algebraic group $G$, which are also all the $\bK$-points of $G$ as $\bK$ is algebraically closed. We denote by $\RGzd$ the resulting smashing colocalization and $\QCohzd(G)$ its image. By definition, $A=\bigoplus_{\lambda\in G(\bK)}\RG_{\lambda}(A)$ for all $A\in\QCohzd(G)$, where $\RG_{\lambda}(A)$ is supported on the closed point $\lambda$. Therefore, $\QCohzd(G)=\prod_{\lambda\in G(\bK)}\QCoh_{\lambda}(G)$. 
\edefn

\bprop[symmetric monoidality of decomposition]\label{prop: symmetric-monoidality-decomposition}
    $\QCohzd(G)\subseteq\QCoh(G)$ inherits a symmetric monoidal structure from the convolution. Moreover, if we identify $\QCoh_{\lambda}(G)$ with $\QCoh_0(G)$ using translation by $\lambda$, then under the identification 
    $$\QCohzd(G)=\prod_{\lambda\in G(\bK)}\QCoh_{\lambda}(G)=\Fun(G(\bK),\QCoh_0(G)),$$ 
    the symmetric monoidal structure corresponds to the Day convolution on the right, namely
    $(A\otimes B)_\lambda=\bigoplus_{\mu+\nu=\lambda}A_\mu\otimes B_\nu$ for $A=(A_\lambda),B=(B_\lambda)\in\QCoh_0(G)^{G(\bK)}$. 
\eprop
\bpf
Let $a\colon G\times G\to G$ be the addition map and let $t_\lambda$ denote translation by $\lambda$. Recall that the convolution product is $A\otimes B=a_*(A\boxtimes B)$. If $A$ is supported at $\mu$ and $B$ at $\nu$, then $A\boxtimes B$ is supported at $(\mu,\nu)$. Its restriction to $a^{-1}(G\setminus\{\mu+\nu\})$ is zero, so $A\otimes B$ is supported at $\mu+\nu$. Since the unit is supported at $0$ and convolution preserves direct sums in each variable, $\QCohzd(G)$ inherits the convolution product.

For $A_\mu,B_\nu\in\QCoh_0(G)$, the identity $a\circ(t_\mu\times t_\nu)=t_{\mu+\nu}\circ a$ gives
\[
(t_\mu)_*A_\mu\otimes(t_\nu)_*B_\nu=(t_{\mu+\nu})_*(A_\mu\otimes B_\nu).
\]
Taking the direct sum over $\mu+\nu=\lambda$ and translating back to $0$ gives the stated formula for $(A\otimes B)_\lambda$. The same identity for the sum maps $G^k\to G$ is compatible with composition and permutations, including the empty sum at $0$. Thus these identifications give the symmetric monoidal equivalence with Day convolution.
\epf
The above proposition also applies to periodic complexes: 

\brem\label{rem: periodic-monoidal-decomposition}
    Recall from Remark \ref{rem: periodic-complexes-as-dg-modules} the algebra $\bK[\beta^\pm]\in\De(\bK)$, viewed as an algebra in $\QCoh_0(G)$ along $0_*$, the pushforward along the unit of $G$. Taking modules in the identification in Proposition \ref{prop: symmetric-monoidality-decomposition}, we get the equivalence
    $$\Mod_{\bK[\beta^\pm]}(\QCohzd(G))=\Fun(G(\bK),\Mod_{\bK[\beta^\pm]}(\QCoh_0(G))),$$
    symmetric monoidally if $|\beta|$ is even or $2=0$ in $\bK$, where we use the Day convolution for the right hand side. This gives the analog of Proposition \ref{prop: symmetric-monoidality-decomposition} for periodic complexes. 
\erem

\brem[functoriality in $G$]
    For a map of algebraic groups $f\colon G\to H$ over $\bK$, it is clear that $f_*$ sends $\QCohzd(G)$ into $\QCohzd(H)$, and specifically for $\lambda\in G(\bK)$ sends the direct summand $\QCoh_\lambda(G)$ to $\QCoh_{f(\lambda)}(H)$. 
\erem
\bpf[Proof of Theorem \ref{thm: main eigenspace decomposition}]
Take $G=\Ga^r=\Spec(A_r)$, where $A_r=\bK[T_1,\ldots,T_r]$. For $\ov{\sX}\in\End^r(\De)$, the underlying homology is finite dimensional in total in the non-periodic case, and over one period in the periodic case. Hence it is supported on a finite set $S$ of closed points. Since $S$ is discrete, the local cohomology counits give $\bigoplus_{\vec{\la}\in S}\RG_{\fm_{\vec{\la}}}(\ov{\sX})\simeq\ov{\sX}$ by \cite[\texttt{0AEF}]{stacks}. Each underlying summand belongs to $\De$ as a retract, and $\De$ is closed under finite sums. Proposition \ref{prop: symmetric-monoidality-decomposition} and Remark \ref{rem: periodic-monoidal-decomposition} therefore restrict to the asserted symmetric monoidal equivalence.
\epf

\section{Sutured instanton homology}\label{sec: Link spectral sequences for sutured instanton homology}
In this section, we review the construction of sutured instanton homology and prove the main theorems.
\subsection{Review of the constructions}\label{eq: review of the constructions}
In this subsection, we review the definition of sutured instanton homology $SHI$ and its surface grading.
\bdefn[{\cite[Definition 2.2]{juhasz2006holomorphic}, \cite[Definition 3.5]{juhasz2008floer} }]\label{defn: balanced sutured manifold}
A \emph{balanced sutured manifold} $(M,\ga)$ consists of a possibly disconnected, oriented, compact 3-manifold $M$ with a non-empty boundary together with a closed oriented 1-submanifold $\ga$ in $\partial{M}$. Let $A(\ga)=[-1,1]\times\ga$ be an annular neighborhood of $\ga\subset \partial{M}$ and let $R(\ga)=\partial{M}\backslash{\rm int}(A(\ga))$. We require the following properties.
\benum
    \item Neither $M$ nor $R(\ga)$ has a closed component.
    \item If $\partial{A(\ga)}=-\partial{R(\ga)}$ is oriented in the same way as $\ga$, then we require that this orientation of $\partial{R(\ga)}$ induces the orientation on $R(\ga)$, which is called the \emph{canonical orientation}.
    \item Let $R_+(\ga)$ be the part of $R(\ga)$ for which the canonical orientation coincides with the induced orientation on $\partial{M}$ from $M$, and let $R_-(\ga)=R(\ga)\backslash R_+(\ga)$. We require that $\chi(R_+(\ga))=\chi(R_-(\ga))$, which is called the \emph{balanced condition}. 
\eenu
Furthermore, if $\chi(R_+(\ga)\cap F)=\chi(R_-(\ga)\cap F)$ for every component $F$ of $\pa M$, then $(M,\ga)$ is called \emph{strongly balanced}.
\edefn

Following Kronheimer--Mrowka \cite[\S 7]{kronheimer2010knots} and Baldwin--Sivek \cite[\S 9]{BS2015naturality}, we define sutured instanton homology as follows.

\bdefn[{\cite[Definitions 9.1 and 9.2]{BS2015naturality}}]\label{defn: closure}
Suppose $(M,\ga)$ is a balanced sutured manifold as in Definition \ref{defn: balanced sutured manifold}. Let $\sD=(Y,R,\omega)$ be a \emph{marked odd closure} (also called a \emph{closure} for short) of $(M,\ga)$, which consists of
\benum
    \item a closed, oriented, connected $3$-manifold $Y$ containing $M$ as a submanifold;
    \item a closed, oriented, connected, embedded surface $R$ of genus $g\ge 2$, called the \emph{distinguished surface}, where $g$ is also called the \emph{genus} of the closure;
    \item $\omega=\al\cup\eta$ for an oriented, connected, embedded curve $\al$ in $Y$ which intersects $R$ transversely once, and a homologically essential, embedded curve $\eta$ in $R$, and hence also in $Y$; we will call $\omega$ the \emph{bundle set}.
\eenu
Here we omit some conditions on the embeddings for simplicity. An \emph{odd closure} $\sD=(Y,R,\omega=\al)$ consists of the data above except for $\eta$. We also call $\al$ the \emph{bundle set}. 
\edefn
\bdefn\label{defn: SHI of closure}
For each odd closure or marked odd closure $\sD$ of $(M,\ga)$, we define \[SHI(\sD)=I(Y,\omega)_{(2,2g(R)-2)}\] to be the simultaneous generalized eigenspace of the $\mu$-operators $(\mu(\pt),\mu(R))$ on instanton Floer homology $I(Y,\omega)$ with eigenvalues $(2,2g(R)-2)$ (see \cite[\S 8]{Ye2026infty} and \cite{donaldson2002floer,DK1992instanton} for the definition of the $\mu$-operators). 
\edefn
\begin{notation}\label{not: SHI}
    Following \cite[Theorem 1.18 and Definition 9.18]{BS2015naturality}, we denote by $\SHI(M,\ga)$ the projectively transitive system of $\mathbb{Z}/2$-graded $\mathbb{C}$-modules $SHI(\sD)$ associated to all marked odd closures $\sD$, which is referred to as the \emph{twisted version of sutured instanton homology}. Note that there is no canonical absolute $\mathbb{Z}/2$-grading; only a relative $\mathbb{Z}/2$-grading is defined, as discussed in \cite[\S 2.6]{kronheimer2010instanton}.

Typically, the notation $SHI(M,\ga)$ refers to the isomorphism class of the $\mathbb{Z}/2$-graded $\mathbb{C}$-modules $SHI(\sD)$ over all odd closures and marked odd closures $\sD$. Although the associated bundle sets may differ, the corresponding modules are still isomorphic. In this paper, however, whenever we write $SHI(M,\ga)$, we fix a marked odd closure and suppress it from the notation. Thus $SHI(M,\ga)$ denotes a fixed $\C$-module. We also fix the auxiliary data used to define cobordism maps, including those associated with families of cobordisms, so that linear maps between such modules are well-defined rather than defined only up to multiplication by units in $\mathbb{C}$. Moreover, $SHI(M,\ga)$ carries a $\mathbb{Z}/2$-grading.

Finally, there is an isomorphism preserving the relative $\mathbb{Z}/2$-gradings,
\[
\SHI(M,\ga)\cong SHI(M,\ga).
\]
Furthermore, all equalities and maps involving $\SHI$ induce corresponding ones involving $SHI$ that are well-defined up to multiplication by units in $\mathbb{C}$.
\end{notation}
Then we focus on the surface grading on $\SHI(M,\ga)$, which was first introduced by Kronheimer--Mrowka in \cite[\S 7]{kronheimer2010knots} and \cite[\S 2.5]{kronheimer2010instanton}. The independence of the auxiliary choices was proved by Li \cite{li2019direct} and Ghosh--Li \cite{li2019decomposition}. See also Wang \cite[\S 4.2.4]{wang2020cosmetic} for alternative proofs.

\bdefn[{\cite[Definition 2.26]{li2019decomposition}}]\label{defn: admissible surface}
Suppose $(M,\ga)$ is a balanced sutured manifold. A surface $S\subset (M,\ga)$ is called \emph{(strongly) admissible} if \benum
    \item it is a properly embedded, oriented surface in generic position;
    \item it may be disconnected and may have closed components other than spheres;
    \item each component of $\pa S$ intersects $\ga$ nontrivially and transversely.
    \item The number \begin{equation}\label{eq: c(S)}
        c(S)=\frac{1}{2}|\pa S\cap \ga|-\chi(S)
    \end{equation}is even, where $|\pa S\cap \ga|$ is the geometric intersection number.
\eenu 
If the last condition is dropped, then we say $S$ is \emph{weakly admissible}. Note that all components of a weakly admissible surface are also weakly admissible.
\edefn
\brem\label{rem: admissible}
We exclude sphere components because many results in \cite[\S 4.2.4]{wang2020cosmetic} require this assumption. We drop the last condition for weakly admissible surfaces because it is a technical assumption in the construction of the grading in sutured instanton theory. In particular, it is a necessary condition for extending the surface $S$ to a closure of $(M,\ga)$. We will also construct gradings for weakly admissible surfaces by stabilizations and renormalizations.
\erem

% In this subsection, we review properties analogous to Lemma \ref{lem: grading properties SFC} for surface gradings in sutured instanton theory. We also review other basic properties about $SHI$. All results in this subsection do not depend on the rest of the paper so they may be used in previous sections.

Recall from Notation \ref{not: SHI} that $\SHI$ is a projectively transitive system of $\C$-modules and $SHI$ is an actual $\C$-module when we fix the closure. Then we collect the following properties from the references above. To relate the grading more directly to the eigenvalues of $\mu$-operators, we multiply the original grading by $2$.

\begin{figure}[htbp]
\centering
\resizebox{0.35\textwidth}{!}{%
\begin{circuitikz}
\tikzstyle{every node}=[font=\LARGE]
\draw [->, >=Stealth] (6.5,12.5) -- (9.75,12.5);
\draw [->, >=Stealth] (9.75,10.75) -- (6.5,10.75);
\draw [->, >=Stealth, dashed] (8,12.25) -- (8,11);
\draw [->, >=Stealth] (4.5,8.75) -- (7.5,8.75);
\draw [->, >=Stealth] (7.5,7) -- (4.5,7);
\draw [->, >=Stealth] (8.75,8.75) -- (12,8.75);
\draw [->, >=Stealth] (12,7) -- (8.75,7);
\draw [->, >=Stealth, dashed] (10.75,9.5) -- (11.5,9.5);
\draw [dashed] (6,8.5) -- (6,8);
\draw [->, >=Stealth, dashed] (5.5,9.25) -- (4.75,9.25);
\draw [dashed] (4.75,9.25) -- (6,7.5);
\draw [dashed] (6,7.5) -- (6,7.25);
\draw [dashed] (10.25,8.5) -- (10.25,8);
\draw [dashed] (10.25,8) -- (10.75,9.5);
\draw [dashed] (11.5,9.5) -- (10.25,7.5);
\draw [dashed] (10.25,7.5) -- (10.25,7.25);
\node [font=\LARGE] at (8.75,11.5) {$\partial S$};
\node [font=\LARGE] at (6,6.25) {$-$};
\node [font=\LARGE] at (10.25,6.25) {$+$};
\node [font=\LARGE] at (9,13) {$\ga$};
\draw [dashed] (6,8) -- (5.5,9.25);
\end{circuitikz}
}%
\caption{The positive ($+$) and negative ($-$) stabilizations of the surface $S$ near the dashed boundary $\partial S$.\label{fig:stabilization}}
\end{figure}
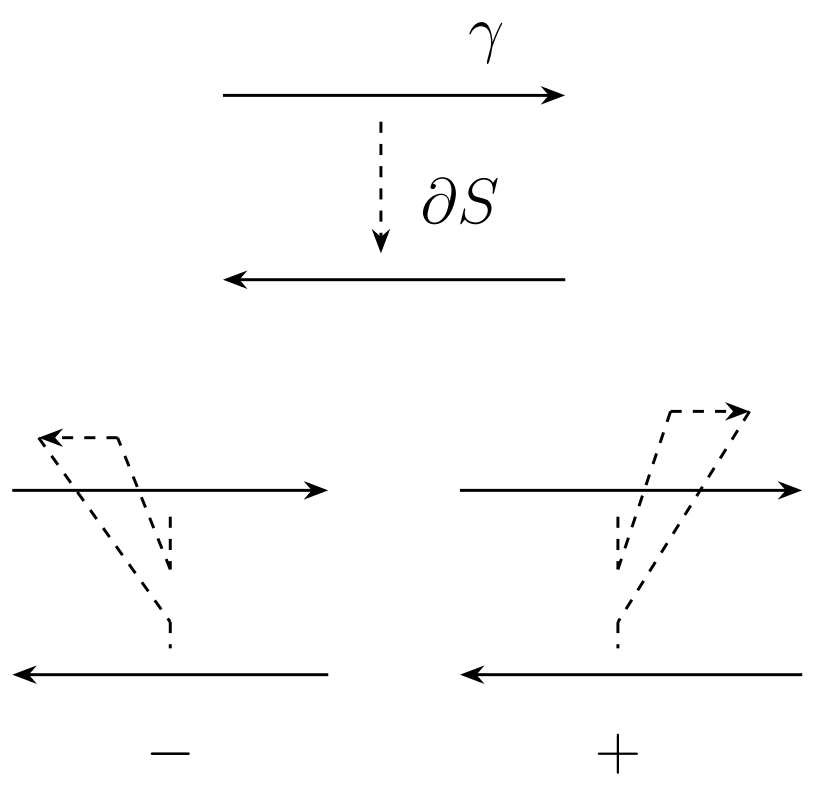

\bprop[{\cite[Theorems 2.21 and 2.27]{LY2020} and \cite[Theorem 2.4]{LY2021enhanced}}]\label{prop: grading properties SHI}
Suppose $(M,\ga)$ is a balanced sutured manifold, and suppose $S$ is a weakly admissible surface as in Definition \ref{defn: admissible surface}. Then there exists a grading on $\SHI(M,\ga)$ indexed by $\la\in\Z$, called the surface grading, satisfying the following properties. \benum
    \item\label{g1} The grading is supported in $c(S)+2\Z$, where $c(S)$ is the quantity in \eqref{eq: c(S)}. In particular, the grading associated to an admissible surface is supported in $2\Z$. Multiple weakly admissible surfaces induce a multi-grading.
    \item\label{g2} A minus sign denotes orientation reversal. Then we have\[\SHI(M,\ga,-S,\la)=\SHI(M,\ga,S,-\la),\]\[\SHI(M,-\ga,S,\la)\cong \SHI(M,\ga,S,-\la),\]\[\SHI(-M,\ga,S,\la)\cong \Hom(\SHI(M,\ga,S,\la),\C),\]where the isomorphisms are canonical and homogeneous with respect to the relative $\Z/2$-gradings.
    \item\label{g3} If $(M,\ga)$ is strongly balanced and if $S'$ is another weakly admissible surface with \[[S',\pa S']=[S,\pa S]\in H_2(M,\pa M;\Z),\]then $S$ and $S'$ define the same relative grading. More precisely, we have\[\SHI(M,\ga,S',\la)=\SHI(M,\ga,S,\la+r(S',t)-r(S,t)),\]where $r(S,t)$ is the rotation number from \cite[Definition 3.8]{juhasz2008floer}, which depends on some trivialization $t$ from the strongly balanced condition. We know the difference $r(S',t)-r(S,t)$ is independent of the choice of $t$. For any balanced $(M, \ga)$, suppose $S'$ is obtained from $S$ by $p_+$ positive stabilizations and $p_-$ negative stabilizations (cf.\ \cite[Definition 3.1]{li2019direct} and Figure \ref{fig:stabilization}). Then $S'$ is still weakly admissible and we have  \[\SHI(M,\ga,S',\la)= \SHI(M,\ga,S,\la-p_++p_-).\]
    \item\label{g4} If $S$ has multiple components $\{S_q\}_{q\in\un{p}}$ for $p\in\Z_+$, then we have \[\SHI(M,\ga,S,\la)=\bigoplus_{\substack{\vec{\la}=(\la_1,\ldots,\la_p)\in\Z^p\\\la=\sum_{q=1}^p\la_q}}\SHI(M,\ga,\{S_{q}\}_{q\in\un{p}},\vec{\la}).\]
    \item\label{g5} 
    If $|\la|>c(S)$, then $\SHI(M,\ga,S,\la)=0$.
    \item\label{g6} If $S$ has no closed components and defines a sutured manifold decomposition $(M,\ga)\overset{S}{\rightsquigarrow}(M',\ga')$, then
\[\SHI(M,\ga,S,c(S))\cong\SHI(M',\ga').\]
\eenu
\eprop
% \brem\label{rem: similar property in SFH}
% Later in \S \ref{subsec: grading in SFH}, especially Lemma \ref{lem: grading properties SFC}, we will follow \cite[\S 4.1]{LY2021enhanced} to construct a surface grading on $SFH(M,\ga)$ and collect similar properties as in Lemma \ref{lem: grading properties SHI} from references. Though the properties are similar, the proofs are quite different.
% \erem

\bpf
If $c(S)$ in (\ref{eq: c(S)}) is even, then the grading is constructed by the eigendecomposition of the $\mu$-operator associated to a closed surface $\ov{S}$ extending $S$ in a closure $\sD$ of $(M,\ga)$, where $\la$ is the eigenvalue and lies in $2\Z$ by \cite[Proposition 7.5]{kronheimer2010knots}. The multi-grading is obtained by taking simultaneous generalized eigenspaces of the $\mu$-operators, which commute with each other. The results in \cite{li2019direct,li2019decomposition} show that this grading is independent of auxiliary choices.

If $c(S)$ is odd, then we apply a positive or negative stabilization to $S$ to make it admissible. Then we define the grading by the renormalization in Term (\ref{g3}). Independence of the choice of stabilization follows from \cite[Corollary 4.18]{wang2020cosmetic}.

% Then we verify terms in Lemma \ref{lem: grading properties SFC} by replacing $SFC$ with $\SHI$.

Term (\ref{g1}) follows from the above definition of the grading.

The first identification in Term (\ref{g2}) follows from the fact that $-\ov{S}$ is a closure of $-S$ and $\mu(-\ov{S})=-\mu(\ov{S})$. The two isomorphisms in Term (\ref{g2}) are based on the fact that if $(Y,R,\omega)$ is a closure of $(M,\ga)$, then $(Y,-R,\omega)$ is a closure of $(M,-\ga)$ and $(-Y,R,\omega)$ is a closure of $(-M,\ga)$. The first isomorphism is $\phi^2$ in \cite[Formula (2.1) and Theorem 2.5]{baldwin2019lspace} or one of the four in \cite[Formula (2.2)]{baldwin2022concordance2} (we choose one so that it is canonical), which depends on the relative $\Z/4$ or $\Z/8$-grading. The second isomorphism is induced by the pairing \[I(Y,\omega)\ot I(-Y,\omega)\to \C\]defined by the canonical trace cobordism map.

The special case of Term (\ref{g3}) for stabilizations follows from \cite[Proposition 4.9]{li2019direct} and \cite[Proposition 4.17]{wang2020cosmetic}. The general case of Term (\ref{g3}) is implicit in \cite[Lemma 4.42]{wang2020cosmetic} and the proof of \cite[Proposition 4.12]{wang2020cosmetic}. It suffices to verify that in the assumption of \cite[Lemma 4.42]{wang2020cosmetic} the difference \[\chi(F\cap R_-(\ga))-\chi(F\cap R_+(\ga))\]coincides with $r(S,t)-r(T,t)$, which can be checked directly from the definition of the rotation number in \cite[Definition 3.8]{juhasz2008floer}.

If every component $S_q$ is admissible, choose extensions $\ov{S}_q$ in a common closure such that their union is an extension $\ov{S}$ of $S$. Then Term (\ref{g4}) follows from
\[\mu(\ov{S})=\sum_{q=1}^p\mu(\ov{S}_q).\]
If some $c(S_q)$ are odd, first stabilize those components and then apply Term (\ref{g3}) to undo the grading shifts.

For admissible $S$, Term (\ref{g5}) follows from \cite[Proposition 7.5]{kronheimer2010knots} and the identity
\begin{equation}\label{eq: c(S) meaning}
-\chi(\ov{S})=c(S).
\end{equation}
If $c(S)$ is odd, apply the admissible case to both a positive and a negative stabilization and use Term (\ref{g3}).

For admissible $S$, Term (\ref{g6}) follows from \cite[Theorem 2.21(2)]{LY2020}. If $c(S)$ is odd, let $S^+$ be a positive stabilization of $S$. The decompositions along $S$ and $S^+$ give diffeomorphic sutured manifolds by \cite[Lemma 2.25(1)]{LY2020}. Since $c(S^+)=c(S)+1$, Term (\ref{g3}) and the admissible case give
\[\SHI(M,\ga,S,c(S))=\SHI(M,\ga,S^+,c(S^+))\cong\SHI(M',\ga').\]
\epf

\subsection{Proofs of the main theorems}
We prove the main theorems from the introduction by combining the results of the previous sections.

\bpf[Proof of Theorem \ref{thm: main grading for cube}]
After stabilizing the $S_l$ near $\pa M$ and away from $L$, choose a common closure $(Y,R,\omega)$ in which they extend to closed surfaces $\ov{S}_l$ (see the proof of \cite[Lemma 5.2]{li2019decomposition}). These stabilizations give the same grading shifts on $M$ and all $M_\ep(L)$ by Proposition \ref{prop: grading properties SHI}\eqref{g3}. Let $L$ denote the induced link in $Y$. Then $(Y_{\ep}(L),R,\omega)$ is a closure of $(M_{\ep}(L),\ga)$ for any $\ep\in\{0,1\}^n$. For $\vec{0}=(0,\dots,0)$, let $L_{\vec{0}}$ denote the link given by the cores of the solid tori in $Y_{\vec{0}}$, framed so that 0-surgery on a core gives the original 1-surgery and 1-surgery gives the original manifold.

We apply \cite[Theorem 8.21]{Ye2026infty} to the surgery cube associated to $L_{\vec{0}}$ with $r=2+k$ and the submanifolds \[\wti{Z}_{1,i}=I\times \pt,~\wti{Z}_{2,i}=I\times R,\aand \wti{Z}_{2+l,i}=I\times \ov{S}_l\text{ for }i\in[1,n].\]Let \[\ov{\sX}\in\End^{r}(\Fun(\bbox^n,\De^c_2(\C)))\] be the resulting hypercube equipped with homotopy-coherently commutative operators. For any $\vec{\la}=(\la_1,\ldots,\la_r)\in\C^r$, let
\[
\sX_{\vec{\la}}=\Eig^{n,r}_{\vec{\la}}(\ov{\sX}):\bbox^n\to\De^c_2(\C)
\]
denote the generalized eigenspace from Definition \ref{defn: eigen functor} associated to the eigenvalues $\vec{\la}$. By Definition \ref{defn: SHI of closure} and Proposition \ref{prop: grading properties SHI}, we consider
\[
\sX_{(2,2g(R)-2)}=\bigoplus_{\vec{\la}'\in\C^k}\sX_{(2,2g(R)-2,\vec{\la}')}.
\]
The spectral sequence then follows from \cite[\S 1.2.2]{HAlurie2017} and the cube filtration on the hypercube (cf.\ \cite[Definition 7.5]{Ye2026infty}). The convergence of the spectral sequence follows from \cite[Theorem 8.21 (3)]{Ye2026infty} together with the iterated procedure described in \cite[Remark 7.18]{Ye2026infty}. Note that from Proposition \ref{prop: grading properties SHI} \eqref{g1}, the admissibility of the surfaces implies that only the summands with $\vec{\lambda}'\in (2\Z)^k$ can be nonzero on homology, giving the surface multi-grading.
\epf

\bpf[Proof of Theorem \ref{thm: main bypass cube}]
From \cite[\S 4]{BS2022khovanov}, after attaching contact $1$-handles to $(M,\ga_\ep)$, the bypass exact triangle can be realized as a special case of the surgery exact triangle on a suitable closure of $(-M,-\ga)$, and the bypass map is the restriction of the corresponding cobordism map. We discuss the choice of bundle sets below.

Note that in the proof of the bypass exact triangle in \cite[\S 4, (16)-(19)]{BS2022khovanov} (see also \cite[\S 8]{LY2025dimension}), one chooses \emph{special} closures of the balanced sutured manifolds at the three vertices so that the original framed knot in the surgery exact triangle bounds a once-punctured torus, with the framing induced by the torus. The existence of such a torus is used to cancel the extra bundle sets on the $3$-manifolds and $4$-dimensional cobordisms introduced in the surgery exact triangle. Although the bundle sets on all three $3$-manifolds can be canceled simultaneously via a homological trick as in \cite[\S 2.2]{baldwin2020concordance}, this can be done for at most two of the three cobordisms.

The complete proof uses the naturality result \cite{BS2015naturality} to relate different choices of closures. Since we only consider the bypass maps between $SHI(M,\ga_\ep)$, we may allow extra bundle sets in one of the other two maps in the bypass exact triangle, so that the bypass exact triangle becomes an actual surgery exact triangle for the chosen special closure. This modification is necessary for our application to the link spectral sequence.

As in the knot case illustrated in \cite[Figure 7]{BS2022khovanov}, we choose a special closure $(-Y_{\vec{1}},-R,\omega)$ of $(-M,-\ga_{\vec{1}})$ in which each component of the surgery link $L_{\vec{1}}\subset -Y_{\vec{1}}$ bounds a punctured torus disjoint from the other components and has the framing induced by that torus. In applying the proof of Theorem \ref{thm: main grading for cube}, we choose the pairing arcs for the stabilized surfaces together with the auxiliary arcs for the punctured tori, so that the same closure realizes the surface gradings. The surgery edges give the bypass maps with the bundle choices above. Each extra meridian is null-homologous after filling, since the corresponding longitude bounds a punctured torus and the $1$-surgery kills the sum of the longitude and meridian. By \cite[Remark 7.19]{Ye2026infty}, we can therefore remove the extra bundle set without changing the surface gradings.

Finally, we dualize the filtered complex and reverse the cube indexing by $\ep\mapsto\vec{1}-\ep$. The dual bypass maps give the first differential. By Term (\ref{g2}) of Proposition \ref{prop: grading properties SHI}, the summand in grading $\vec{\la}$ corresponds to the summand in grading $-\vec{\la}$.
\epf
\brem
In the proof of Theorem \ref{thm: main bypass cube}, if it is not necessary for the differential to be given by bypass maps, we can alternatively use \emph{any} closure of $(-M,-\ga)$ with sufficiently large genus $g(R)$ to carry out the proof, at the cost of allowing extra bundle sets on both the manifolds and cobordisms. In this case, the restrictions of the cobordism maps are no longer bypass maps in the sense of \cite[\S 4]{BS2022khovanov}.
\erem

\section{Mixed hypercube}\label{sec: mixed hypercube}
In this section, we use the special closure from Assumption \ref{assump: topological} to prove Theorem \ref{thm: mixed hypercube}. We adopt the notation in \S \ref{subsec: combination of bypass maps} and \S \ref{eq: review of the constructions}, especially that in Assumption \ref{assump: topological}. Note that a closure $(Y_{\vec{0}},R)$ of $(M,\ga_{\vec{0}})$ is indeed a marked odd closure $(Y_{\vec{0}},R,\omega)$ from Definition \ref{defn: closure} with the bundle set $\omega$ as additional data. Since $\omega$ is $1$-dimensional, it can be made disjoint from all surgery curves $C_i^\pm$ associated to the bypass maps. 

% For the closed $3$-manifold $Y^s_\ep$ from Assumption \ref{assump: topological}, we write $\sD^s_\ep=(Y^s_\ep,R,\omega)$ for the induced marked odd closure of $(M,\ga)$. 

Following \cite[Construction 7.10]{Ye2026infty}, we will write $\bY=(Y,\omega)$, $\bW=(W,\nu_W)$, and $\bX=(X,\nu_X)$ for manifolds and cobordisms with bundle sets. For simplicity, we also write $\un{n}=\{1,\dots,n\}$. 

For the composition law, the subgroup $\Phi$ on $W=W_1\circ W_0$ is the largest subgroup of $H^1(W;\Z/2)$ whose restriction to each $W_i$ is $\Phi_i$, as in the introduction of \cite{Ye2026infty}. In our case, each connected component of $W_0$ has exactly one outgoing boundary component, since surgery preserves components and excision only merges them. Thus $H^0(W_0;\Z/2)\to H^0(Y;\Z/2)$ is surjective for the common boundary $Y$. The Mayer--Vietoris sequence with $\Z/2$ coefficients gives
\[
H^0(W_0)\oplus H^0(W_1)\longrightarrow H^0(Y)\xrightarrow{\delta}H^1(W)\xrightarrow{\rho}H^1(W_0)\oplus H^1(W_1).
\]
The first map is surjective, so $\delta=0$ and $\rho$ is injective. If $\Phi_0=\Phi_1=0$, the largest subgroup with these restrictions is $\ker\rho=0$. This applies at each stage of the composition, so we may take all cobordism subgroups $\Phi$ and their boundary restrictions $\phi$ to be trivial.

\subsection{Step I: Hyperbox complexes}\label{subsec: box complex}

In this step, we construct $2^n$ hypercubes for all $s\in\{+,-\}^n$ for closed $3$-manifolds.

For $s=\vec{+}=(+,\dots,+)$, by \cite[Construction 7.11]{Ye2026infty}, we obtain a hypercube \begin{equation}\label{eq: surgery cube plus}
    \sX^{\vec{+}}=\sX[\bY_{\vec{0}}\sqcup \bigsqcup_{j=1}^n \bY_{\varphi_j},\{C_i^+\}_{i\in\un{n}}]:\bbox^n\to\BI^\pa,
\end{equation} where $\BI^\pa$ denotes the $\infty$-category of instanton cobordisms with full boundary conditions introduced in \cite[Definition 5.6]{Ye2026infty}. Moreover, composing $\sX^{\vec{+}}$ with the $\infty$-categorical instanton functor \[CI:\BI^\pa\to\De_2^c(\Z)\] from \cite[Remark 6.2]{Ye2026infty} and then applying the derived tensor functor $-\otimes_{\Z}^{\mathbf L}\C$ gives a hypercube complex in $\De_2^c(\C)$. For simplicity, we write $\De=\De_2^c(\C)$ throughout this section.

For general $s\in\{+,-\}^n$, Assumption \ref{assump: topological} \eqref{as4} and \cite[Construction 7.11]{Ye2026infty} determine a hypercube
\begin{equation}\label{eq: surgery cube general s}
   \sX^s=\sX[\bY_{\vec{0}}\sqcup \bigsqcup_{j=1}^n \bY_{\varphi_j},\{\bW^{s_i,\prime}_i\}_{i\in\un{n}}]:\bbox^n\to\BI^\pa. 
\end{equation}
The above hypercube is the compression of the hyperbox from \cite[Remarks 2.3 and 7.13]{Ye2026infty}. Also note that, except for $\vec{0}\in\bbox^n$, the images of the other objects in $\bbox^n$ may vary with $s$, according to which mapping tori
$\bY_{\varphi_j}$ remain as separate components. Again, we apply $CI$ to obtain a hypercube complex. For simplicity, we abuse notation and continue to use $\sX^s$ to denote the resulting hypercube complex.

\subsection{Step II: Sutured instanton homology}\label{subsec: sutured instanton homology}

In this step, we work over $\C$ and pass to sutured instanton homology as in Definition \ref{defn: SHI of closure} and the proof of Theorem \ref{thm: main grading for cube}, so that each mapping torus contributes a 1-dimensional top eigenspace. By fixing homogeneous generators, we obtain $2^n$ hypercube complexes before passing to the surface multi-grading.

For $\vec{+}\in\{+,-\}^n$, Assumption \ref{assump: topological} \eqref{as2} ensures that the surgery curves $\{C_i^+\}_{i\in\un{n}}$ associated with the positive bypass maps are disjoint from the embedded surface $R$. We also suppose $\pt$ is a point in $Y_{\vec{0}}$ that is disjoint from all surgery curves and the surface $R$. We also use $\pt$ and $R$ to denote the induced submanifolds in the closures after surgery. As in Assumption \ref{assump: topological} \eqref{as3}, write $R_{\varphi_i}$ for the fiber surface in the mapping torus $Y_{\varphi_i}$ and for the induced surface after the canonical cobordism. Similarly, we write $\pt_{\varphi_j}$ for a point in $Y_{\varphi_j}$ that is disjoint from $R_{\varphi_j}$ and the region of the canonical cobordism, and for the induced point after this cobordism.

As in the proof of Theorem \ref{thm: main grading for cube}, we apply \cite[Theorem 8.21]{Ye2026infty} to the surgery cube \begin{equation}\label{eq: eq surgery cube plus 2}
    \sX[\bY_{\vec{0}}\sqcup \bigsqcup_{j=1}^n\bY_{\varphi_j},\{C_i^+\}_{i\in\un{n}}]
\end{equation}in \eqref{eq: surgery cube plus}, with $r=2+2n$ and the submanifolds \begin{equation}\label{eq: cob of surfaces 1}
    \wti{Z}_{1,i}=I\times \pt, ~\wti{Z}_{2,i}=I\times R, ~\wti{Z}_{2j+1,i}=I\times \pt_{\varphi_j},~\wti{Z}_{2j+2,i}=I\times R_{\varphi_j}\text{ for }i,j\in\un{n}.
\end{equation} Then we obtain homotopy-coherently commutative operators on $\sX^{\vec{+}}$\begin{equation}\label{eq: lift}
    \ov{\sX}^{\vec{+}}\in\End^{2+2n}(\Fun(\bbox^n,\De)).
\end{equation}Following Definition \ref{defn: eigen functor}, for $g=g(R)=g(R_{\varphi_j})$ and $\vec{\la}_0=(2,2g-2,\dots,2,2g-2)$, let \begin{equation}\label{eq: eigen plus}
    \sX^{\vec{+}}_{\vec{\la}_0}: \bbox^n\to \De
\end{equation}denote the corresponding generalized eigenspace. From \cite[Theorem 7.7 and Proposition 7.4]{kronheimer2010knots}, excision on mapping tori with monodromy $\varphi_j$ and $\varphi_j^{-1}$ implies\begin{equation}\label{eq: dim 1}
    I(\bY_{\varphi_j})_{2,2g-2}\cong \C.
\end{equation}Fix a homogeneous generator of $I(\bY_{\varphi_j})_{2,2g-2}$ for each $j$ and use the resulting isomorphisms to remove the mapping torus factors on every subcube. In particular, for $s=\vec{+}$ this identifies $\sX^{\vec{+}}_{\vec{\la}_0}$ with the $(2,2g-2)$-generalized eigenspace of \[\sX[\bY_{\vec{0}},\{C_i^+\}_{i\in\un{n}}].\]

The case of general $s\in\{+,-\}^n$ is similar. By composing with $CI$ and then taking the tensor product with $\C$, we regard the hypercube $\sX^s$ in \eqref{eq: surgery cube general s} as
\[
\sX^s:\bbox^n\to\De^c_2(\C).
\]
By applying \cite[Theorem 8.21]{Ye2026infty} to \eqref{eq: cob of surfaces 1}, we obtain $2+2n$ homotopy-coherently commutative operators on $\sX^s$. Upon passing to the generalized eigenspace associated with $\vec{\la}_0=(2,2g-2,\dots,2,2g-2)$, it is important to observe that $R_{\varphi_j}$ is isotopic to $R$ in the target of the canonical cobordism in \eqref{eq: canonical cobordism 2}, and similarly $\pt_{\varphi_j}$ is isotopic to $\pt$. By Corollary \ref{cor: nilpotent}, the $(2,2g-2,2,2g-2)$-generalized eigenspace of the operators associated with $(\pt,R,\pt_{\varphi_j},R_{\varphi_j})$ coincides with the $(2,2g-2)$-generalized eigenspace of the first two operators. 

Using the fixed homogeneous generators in \eqref{eq: dim 1}, we remove the remaining mapping torus factors and obtain the hypercube complexes
\begin{equation}\label{eq: eigenspaces for all s}
    \{\sX^s_{\vec{\la}_0}\}_{s\in\{+,-\}^n}:\bbox^n\to\De.
\end{equation}
With the same auxiliary data on common subcubes, these complexes agree on those subcubes, so the complexes satisfy Assumption \ref{assump: mixed condition} before imposing the surface multi-grading. Note that the operators also agree on those subcubes.
\subsection{Step III: Surface multi-grading}\label{subsec: gradings}
In the final step, we use the cobordisms of surfaces \eqref{eq: Z for S general} associated to admissible surfaces to further decompose the hypercube complexes \eqref{eq: eigenspaces for all s} and show that the further decompositions still satisfy the mixed condition, which is exactly the statement in Assumption \ref{assump: mixed condition}.

As in \eqref{eq: lift}, for $s=\vec{+}$, we apply \cite[Theorem 8.21]{Ye2026infty} to the surgery cube \eqref{eq: surgery cube plus} with $r=2+2n+k$ and the submanifolds in \eqref{eq: cob of surfaces 1} and \begin{equation}\label{eq: choice of Z grading plus}
  \wti{Z}_{2+2n+l,i}^{\vec{+}}=I\times \ov{S}^{\vec{+}}_l\text{ for } i\in\un{n}, l\in\un{k}.
\end{equation}from Assumption \ref{assump: topological} \eqref{as1} and \eqref{eq: Z for S general}. Then there are further homotopy-coherently commutative operators beyond (\ref{eq: lift}), written as\begin{equation}\label{eq: further lift}
    \ovv{\sX}^{\vec{+}}\in\End^{2+2n+k}(\Fun(\bbox^n,\De)).
\end{equation}

Then we obtain a further eigendecomposition of \eqref{eq: eigen plus}, written as \begin{equation}\label{eq: decomp X sut}
\sX^{\vec{+}}_{\vec{\la}_0}=\bigoplus_{\vec{\la}\in\C^k}\sX^{\vec{+}}_{(\vec{\la}_0,\vec{\la})}.
\end{equation}
From Proposition \ref{prop: grading properties SHI} \eqref{g1}\eqref{g5}, we know that the direct summand is nontrivial only for finitely many $\vec{\la}\in(2\Z)^k$. 

The cobordisms of surfaces for general $s\in\{+,-\}^n$ are no longer product cobordisms, but involve those from \eqref{eq: Z for S general}. The existence of the homotopy-coherently commutative operators is ensured by Assumption \ref{assump: topological} \eqref{as4}. Note that here we only need the fact that \begin{equation}
    \label{eq: all coborodisms weak}
        \big\{(W^{s_i,\prime}_i,\{\wti{Z}^{s,\prime}_{i,l}\}_{l\in[1,k]})\big\}_{i\in[1,n]}
    \end{equation}are supported in disjoint regions, which is weaker than \eqref{eq: all coborodisms}. The full condition in \eqref{eq: all coborodisms} will be used to verify the mixed condition.

For $i\in\un{n}$ and $l\in\un{k}$, we take
\begin{equation}\label{eq: choice of Z grading}
    \wti{Z}_{2+2n+l,i}^s=\wti{Z}^{s,\prime}_{i,l}
\end{equation}
from \eqref{eq: composition of cobordism surfaces}, with product cylinders over the surfaces in the mapping torus components unused in direction $i$. The resulting hypercube complex with operators is written as
\begin{equation}\label{eq: further lift 2}
    \ovv{\sX}^{s}\in \End^{2+2n+k}(\Fun(\bbox^n,\De)).
\end{equation}
As in \eqref{eq: decomp X sut} and the discussion after it, we have
\begin{equation}\label{eq: decomp X sut general s}
\sX^{s}_{\vec{\la}_0}=\bigoplus_{\vec{\la}\in(2\Z)^k}\sX^{s}_{(\vec{\la}_0,\vec{\la})}
\end{equation}and only finitely many summands are nontrivial. 

Moreover, since $F_{i,l}^s$ is a collection of tori in a mapping torus, from \cite[Proposition 7.5]{kronheimer2010knots}, we know that $0$ is the only eigenvalue of $\mu(F_{i,l}^s)$ on the summand in \eqref{eq: dim 1}. In particular, the generalized eigenspace \eqref{eq: dim 1} is contained in the generalized eigenspace of $\mu(F_{i,l}^s)$ with eigenvalue $0$. Hence the eigendecomposition associated to \[\mu(\ov{S}^{s}_l\sqcup F_{i,l}^s)=\mu(\ov{S}^{s}_l)\ot \id+\id\ot \mu(F_{i,l}^s)\]coincides with the eigendecomposition associated to $\mu(\ov{S}^{s}_l)$. We remove the remaining mapping torus factors using the same homogeneous generators as in Step II.

Now we verify Assumption \ref{assump: mixed condition}. We provide the details because the operators corresponding to different $s\in\{+,-\}^n$ are no longer identical. We first introduce the following definition.

\bdefn
Let $\tbox^n=(0\rightrightarrows 1)^n$, with the two arrows in each coordinate labelled $+$ and $-$. A \emph{signed subcube} $F$ of $\tbox^n$ is obtained by fixing some coordinates at $0$ or $1$ and choosing one arrow in each remaining coordinate. Write $\bbox_F\subset\tbox^n$ for the resulting sub-$\infty$-category and $\cF$ for the poset of signed subcubes ordered by inclusion.
\edefn

We have
\begin{equation}\label{eq: signed face descent}
\tbox^n=\colim_{F\in\cF}\bbox_F,~
\Fun(\tbox^n,\De)=\lim_{F\in\cF^{\opp}}\Fun(\bbox_F,\De).
\end{equation}
For $n=1$, the first identity is the pushout of two copies of $\Delta^1$ along their common boundary. The general case follows because products in the $\infty$-category of $\infty$-categories $\Cat$ preserve colimits separately in each variable. Applying $\Fun(-,\De)$ gives the second identity \cite[\texttt{01FJ}, \texttt{05UR}]{kerodon}.

To construct the mixed cube, we therefore give an object $\sX_F\in\Fun(\bbox_F,\De)$ for each $F\in \cF$ and isomorphisms $\sX_F\to\sX_G|_{\bbox_F}$ for $F\subset G$, compatible along chains of subcube inclusions.

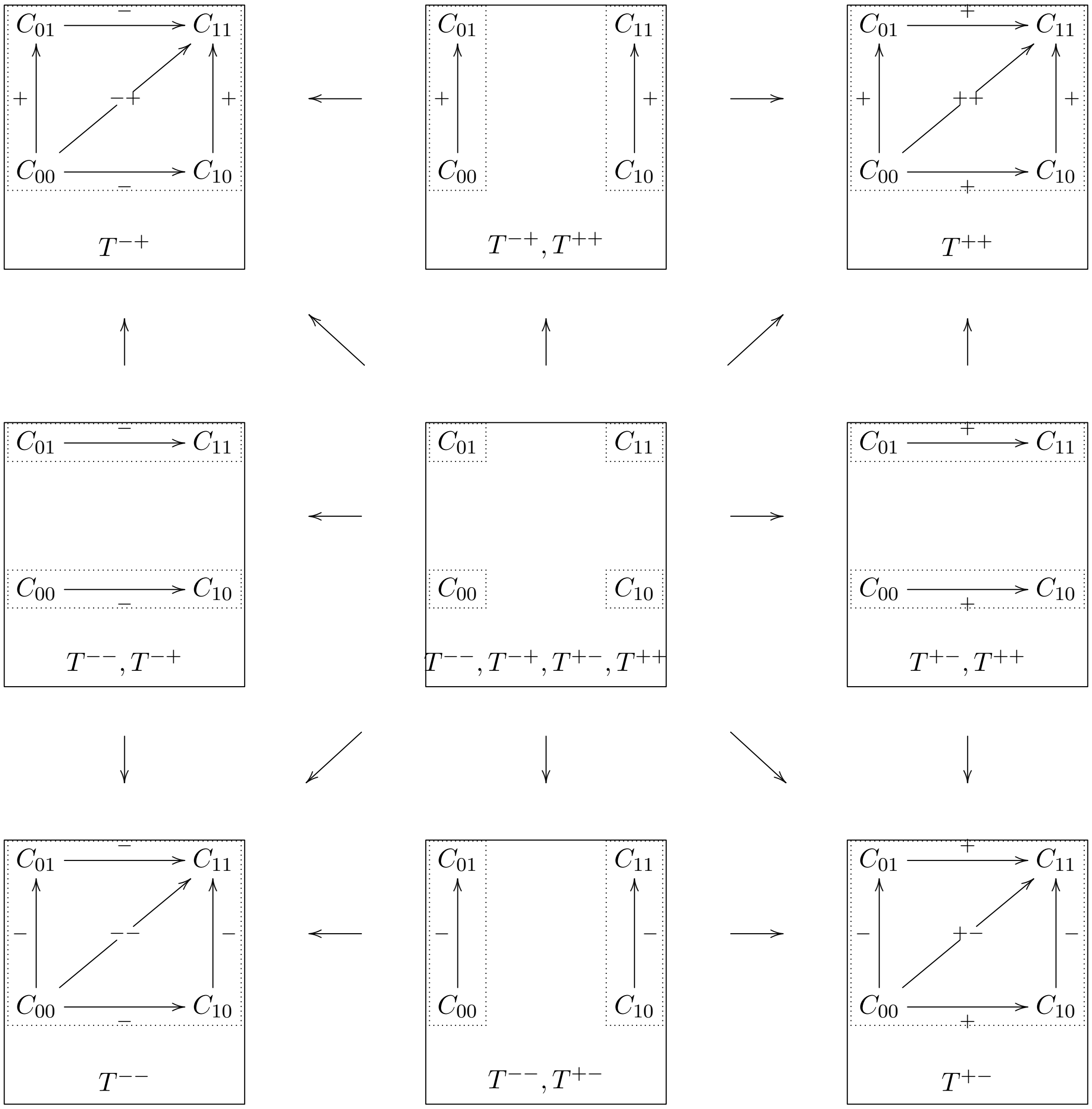
\begin{figure}[htbp]
\centering

\resizebox{0.95\textwidth}{!}{$\xymatrix@C=1.8pc@R=1.6pc{\save [].[dddrr] *+<0.4cm,0.2cm>[F]\frm{}\restore \save [].[rr].[ddrr]*[F.]\frm{}\restore C_{01}\ar[rr]^{-}&&C_{11}&&&\save [].[dddrr] *+<0.4cm,0.2cm>[F]\frm{}\restore \save [].[dd]*[F.]\frm{}\restore C_{01}&&\save [].[dd]*[F.]\frm{}\restore C_{11}&&&\save [].[dddrr] *+<0.4cm,0.2cm>[F]\frm{}\restore  \save [].[rr].[ddrr]*[F.]\frm{}\restore C_{01}\ar[rr]^{+}&&C_{11}\\&&&\ar@{<-}[r]&&&&&&\ar@{<-}[l]&&&\\
    C_{00}\ar[rr]_{-}\ar[uu]^{+}\ar[uurr]|-{-+}&&C_{10}\ar[uu]_{+}&&& C_{00}\ar[uu]^{+}&&C_{10}\ar[uu]_{+}&&&C_{00}\ar[rr]_{+}\ar[uu]^{+}\ar[uurr]|-{++}&&C_{10}\ar[uu]_{+}\\&\save[] *{T^{-+}} \restore &&&&&\save[] *{T^{-+},T^{++}}\restore&&&&&\save[] *{T^{++}}\restore &\\&\ar@{<-}[d]&&\ar@{<-}[dr]&&&\ar@{<-}[d]&&&\ar@{<-}[dl]&&\ar@{<-}[d]&\\&&&&&&&&&&&&\\\save [].[dddrr] *+<0.4cm,0.2cm>[F]\frm{}\restore  \save [].[rr].[rr]*[F.]\frm{}\restore 
     C_{01}\ar[rr]^{-}&&C_{11}&&&\save [].[dddrr]*+<0.4cm,0.2cm>[F]\frm{}\restore  \save [].*[F.]\frm{}\restore C_{01}&&\save [].*[F.]\frm{}\restore C_{11}&&&\save [].[dddrr] *+<0.4cm,0.2cm>[F]\frm{}\restore  \save [].[rr].[rr]*[F.]\frm{}\restore C_{01}\ar[rr]^{+}&&C_{11}\\&&&\ar@{<-}[r]&&&&&&\ar@{<-}[l]&&&\\\save [].[rr].[rr]*[F.]\frm{}\restore C_{00}\ar[rr]_{-}&&C_{10}&&& \save [].*[F.]\frm{}\restore C_{00}&&\save [].*[F.]\frm{}\restore C_{10}&&&\save [].[rr].[rr]*[F.]\frm{}\restore C_{00}\ar[rr]_{+}&&C_{10}\\&\save[] *{T^{--},T^{-+}} \restore &&&&&\save[] *{T^{--},T^{-+},T^{+-},T^{++}}\restore &&&&&\save[] *{T^{+-},T^{++}}\restore &\\&&&&&&&&&&&&\\&\ar@{<-}[u]&&\ar@{<-}[ur]&&&\ar@{<-}[u]&&&\ar@{<-}[ul]&&\ar@{<-}[u]&\\\save [].[dddrr] *+<0.4cm,0.2cm>[F]\frm{}\restore  \save [].[rr].[ddrr]*[F.]\frm{}\restore 
     C_{01}\ar[rr]^{-}&&C_{11}&&&\save [].[dddrr] *+<0.4cm,0.2cm>[F]\frm{}\restore \save [].[dd]*[F.]\frm{}\restore C_{01}&&\save [].[dd]*[F.]\frm{}\restore C_{11}&&&\save [].[dddrr] *+<0.4cm,0.2cm>[F]\frm{}\restore  \save [].[rr].[ddrr]*[F.]\frm{}\restore 
    C_{01}\ar[rr]^{+}&&C_{11}\\&&&\ar@{<-}[r]&&&&&&\ar@{<-}[l]&&&\\
    C_{00}\ar[rr]_{-}\ar[uu]^{-}\ar[uurr]|-{--}&&C_{10}\ar[uu]_{-}&&& C_{00}\ar[uu]^{-}&&C_{10}\ar[uu]_{-}&&&C_{00}\ar[rr]_{+}\ar[uu]^{-}\ar[uurr]|-{+-}&&C_{10}\ar[uu]_{-}\\&\save[] *{T^{--}} \restore &&&&&\save[] *{T^{--},T^{+-}}\restore &&&&&\save[] *{T^{+-}}\restore &}$}

\caption{An example of subcubes and operators when $n=2$ and $k=1$.}
\label{fig:mixed example}

\end{figure}

Consider the case $n=2$ and $k=1$ in Figure \ref{fig:mixed example}. Each vertex lies in all four signed squares, so we take the simultaneous generalized eigendecomposition for all four operators $T^s$ at that vertex. Each signed edge lies in two squares, so we take the simultaneous generalized eigendecomposition for the two corresponding operators on the edge diagram. At each endpoint of a signed edge, each simultaneous generalized eigenspace for all four operators includes naturally into the summand of the edge diagram at that endpoint with the same eigenvalues for the two operators on the edge. These inclusions are compatible along chains of signed subcubes. For a vertex contained in an edge of a signed square, the inclusion into the square summand at that vertex agrees with the composite through the edge summand.

Then we consider the general case. Fix a signed subcube $F$ with initial vertex $\ep$ and terminal vertex $\epp$. For $s\in\{+,-\}^n$ extending the signs of $F$, consider the restriction
\begin{equation}\label{eq: subcubes}
    \sX^s_{\vec{\la}_0}|_{\bbox_F}:\bbox_F\to\De.
\end{equation}
For any two such sign vectors $s,t$, we will show that these restrictions are canonically isomorphic up to a shift in the surface multi-grading. That is, we have
\begin{equation}\label{eq: subcubes with gradings}
    \sX^s_{(\vec{\la}_0,\vec{\la})}|_{\bbox_F}\cong
    \sX^t_{(\vec{\la}_0,\vec{\la}+\vec{\al}^{s,t})}|_{\bbox_F},
\end{equation}
where \eqref{eq: many shift grading} gives
\[
    \vec{\al}^{s,t}=\sum_{\substack{\ep_i=\ep_i^\prime=0\\t_i=-}}f_i-\sum_{\substack{\ep_i=\ep_i^\prime=0\\s_i=-}}f_i.
\]

The first $2+2n$ operators
\begin{equation}\label{eq: operators}
    \mu(\pt),~\mu(R),~\mu(\pt_{\varphi_1}),~\mu(R_{\varphi_1}),~\dots,~\mu(\pt_{\varphi_n}),~\mu(R_{\varphi_n})
\end{equation}
are independent of the choice of $s$ extending the signs of $F$. With the fixed identifications of the mapping torus factors, the restrictions in \eqref{eq: subcubes} agree before imposing the surface multi-grading.

On this common subcube, take the surface operators $T_l^s$ from \eqref{eq: choice of Z grading} for all $s$ extending the signs of $F$. They commute homotopy-coherently by Assumption \ref{assump: topological} (\ref{as4}) and \cite[Theorem 8.21]{Ye2026infty}. We use the same auxiliary data on common subcubes, so for $F\subset G$ the operators on $\bbox_G$ restrict to the corresponding subfamily on $\bbox_F$. Let $\sX_F$ be the direct sum of the simultaneous generalized eigenspaces for the surface operators on $\bbox_F$.

By \eqref{eq: many shift grading}, at each vertex of $\bbox_F$ a simultaneous generalized eigenspace can be nonzero only if $\vec{\la}^{\,t}=\vec{\la}^{\,s}+\vec{\al}^{s,t}$. The shift is independent of the vertex because $s$ and $t$ agree in the varying coordinates of $F$. By Corollary \ref{cor: main eigen}, on forgetting some operators the generalized eigenspace for the remaining operators is the direct sum of the simultaneous generalized eigenspaces with those eigenvalues. The eigenvalues for any one $s$ determine all the others on a nonzero summand. Thus the natural inclusion into the corresponding generalized eigenspace for that $s$ is an isomorphism at each vertex, hence in $\Fun(\bbox_F,\De)$. This gives \eqref{eq: subcubes with gradings}.

For $F\subset G$, forgetting the operators that occur only on $\bbox_F$ gives the same natural inclusions into the generalized eigenspaces on $\bbox_G$ restricted to $\bbox_F$. These inclusions are compatible along chains of subcube inclusions by successive forgetting of operators, as in Figure \ref{fig:mixed example}. Taking direct sums gives the isomorphisms $\sX_F\to\sX_G|_{\bbox_F}$ and hence the mixed cube by \eqref{eq: signed face descent}.

After the regrading in \eqref{eq: many shift grading}, choose a single dg nerve representative of the glued diagram with the graded $SHI$ groups as its vertex complexes. Then each $D^s_{\ep,\epp}$ depends only on the corresponding signed interval (cf.\ \cite[Construction 7.7]{Ye2026infty}), so Assumption \ref{assump: mixed condition} holds.

\bibliographystyle{alpha}

\end{document}